\documentclass[11pt, a4paper]{article}
\usepackage[utf8]{inputenc}
\usepackage[T2A]{fontenc}
\usepackage{mathtext}
\usepackage{mathrsfs}
\usepackage{amssymb}
\usepackage{amsfonts}
\usepackage{amscd}
\usepackage{amsmath, amsthm}
\usepackage{mathtools}
\usepackage{hyperref}
\usepackage[nice]{nicefrac}
\usepackage[symbol]{footmisc}

\DeclareMathOperator*{\res}{Res}

\begin{document}

\begin{center}\textbf{ON POLYNOMIALS OF BINOMIAL TYPE, EXPONENTIAL INTEGRAL AND THE INVERSE LOGARITHMIC DERIVATIVE EIGENPROBLEM}
\end{center}
\begin{center}
           \begin{footnotesize}
	DANIL KROTKOV\footnote[2]{Higher School of Economics, Faculty of Mathematics, e-mail: dikrotkov@edu.hse.ru}
	\end{footnotesize}
\end{center}

\begin{footnotesize}
\noindent\textsc{Abstract}. In this work we continue to study the properties of polynomials of binomial type and their canonical continuations to the complex index by exploring the properties of transformation $\mathfrak{T}\coloneqq1/\mathrm{dlog}$ which acts on formal power series $f(x)$ of the form $x+x^2\mathbb{C}[[x]]$.
\end{footnotesize}

\renewcommand{\contentsname}{Contents}
\tableofcontents

\begin{footnotesize}
\section*{Introduction}

\noindent In the previous paper we considered polynomials of binomial type, i.e. the sequences of \mbox{polynomials} with property $p_n(\alpha+\beta)=\sum_{k=0}^{n}\binom{n}{k}p_k(\alpha)p_{n-k}(\beta)$, which we call \emph{associated to formal power series} \mbox{$f(x)\in x+x^2\mathbb{C}[[x]]$}, if the property $f(\partial/\partial \alpha)p_n(\alpha)=np_{n-1}(\alpha)$ holds (such a series always exists for polynomials of binomial type, see \cite{Dlt}). Also, we studied their canonical continuations to the complex index. Particularly, we were interested in polynomials $\nu_n(\alpha)$, associated to the function $xe^{\mathrm{Ei}(x)-\ln|x|-\gamma}$, and \mbox{having} the property $\alpha\nu'_n(\alpha-1)=n\nu_n(\alpha)$. As it was demonstrated, this property is the direct corollary of the fact, that the function $xe^{-x}$ is the image of the function $xe^{\mathrm{Ei}(x)-\ln|x|-\gamma}$ under the transformation $\mathfrak{T}$:$~f \to f/f'$. Given the potential benefits, the author suggested to look at the sequences of the form $...\to f \to f/f' \to ff'/(f'^2-ff'') \to...$, which we call \emph{$\mathfrak{T}$-chains}.\\
In Section 1 the general question about ``eigenfunctions'' of nonlinear operator $\mathfrak{T}$ is raised, which was solved in a particular case of periodic $\mathfrak{T}$-chains in the previous paper.\\
In Section 2 possible ways to generalize the trivial identity $(\ln\alpha)'=\alpha^{-1}$ in the domain of continuations of binomial polynomials are discussed.\\
In the final third section we provide a different perspective on polynomials of binomial type of integer index, based on the properties of ``derivation'' and ``shift'' operators. We then apply the acquired knowledge to obtain new formulae for the Taylor series of function $(xe^{\mathrm{Ei}(x)-\ln|x|-\gamma})^{inv}$.
\end{footnotesize}

\newpage

\hrule

\section{On the equation $\mathfrak{T}f(x)=\dfrac{f(px)}{p}$}

\hrule

~\\
~\\
In this section we study the equation $f(x)/f'(x)=f(px)/p$ in the domain $x+x^2\mathbb{C}[[x]]$. As it was demonstrated in the previous paper, in case $p=-1$ the solution of this equation takes the form $f(x)=(e^{Ax}-1)/A$ for some constant $A$, since from the condition $p=-1$ it follows that $\mathfrak{T}^2f=f$, which means the $\mathfrak{T}$-chain for such $f$ is periodic, and hence it must be the series of exponent. This section is divided into several parts. At first, we provide the main properties of solutions of this equation. Then the particular cases $p^2=-2$ and $p^4=-4$ are considered. And finally, we discuss general properties of the supplementary series $\varphi_{n}(x)$, which are defined below.

\begin{center}
	\textbf{The main properties of solutions}
\end{center}
\textbf{Proposition 1.1}.
$$
\mathfrak{T}f(x)=\frac{f(px)}{p} ~\Rightarrow~ \exists n \in \mathbb{N}: p^{n}=-n \eqno (1.1)
$$
\emph{Proof}: Suppose such $f(x)$ is nontrivial, i.e. there exists the least $n \geqslant 1$, such that the coefficient of $x^{n+1}$ is not equal to zero. The first $n$ terms of the series $f(x)$ determine the first $n$ terms of the series $\mathfrak{T}f$, thus
$$
f \approx x+ \alpha x^{n+1} ~\Rightarrow~ \mathfrak{T}f \approx x-n\alpha x^{n+1}
$$
But
$$
\mathfrak{T}f \approx x+\alpha p^{n} x^{n+1}, ~ \alpha \neq 0 ~\Rightarrow~ p^{n}=-n
$$\qed\\
\textbf{Proposition 1.2}. Suppose $f \approx x+ \alpha x^{n+1}$. Then it follows that
\begin{align*}
\tag{1.2}\mathfrak{T}f(x)=\frac{f(px)}{p} ~\Rightarrow~ f(x)=x+A_1 x^{n+1}+A_2 x^{2n+1}+...\\
\text{for some}~A_i.
\end{align*}
\emph{Proof}: Suppose $f(x)=x+B_1 x^2+B_2 x^3+...$. Then, since $f'(x)f(px)=pf(x)$, we have the following relation
$$
B_k (k+p^k)=-\sum_{m=1}^{k-1}(m+1)p^{k-m}B_m B_{k-m}
$$
As it is assumed, we have $B_1=B_2=...=B_{n-1}=0$, and $B_n$ is the first nonzero coefficient. Then it follows from the recurrence relation that $B_n (n+p^n)=0$. We put $A_1\coloneqq B_n$. Now notice that since $n+p^n=0$, we have $k+p^k\neq 0 ~~\forall k>n$. That means $B_{n+1}, B_{n+2},...,B_{2n-1}$ are equal to zero as the sums of zero components. So the first nontrivial relation is $B_{2n}(2n+p^{2n})=-B_n^2 p^{n}(n+1)$. $2n+p^{2n}\neq 0$, hence the coefficient $A_2\coloneqq B_{2n}$ is determined by $B_{n}$. Similarly, $B_{2n+1}, B_{2n+2},...,B_{3n-1}$ are equal to zero as the sums of zero components, and the next nontrivial component is $B_n B_{2n}$, which takes place in the expansion of $B_{3n}$. Using the same argument on each step, we obtain that $f(x)$ is of the form $x+A_1 x^{n+1}+A_2 x^{2n+1}+...$ for some $A_i$. \qed
\newpage
\noindent\textbf{Proposition 1.3}.
$$
\mathfrak{T}f(x)=\frac{f(px)}{p},~p^n=-n, ~ f(x)=x+A_1 x^{n+1}+A_2 x^{2n+1}+...
$$ Then the following holds:\\
$$
(k-(-n)^{k-1})A_k = \sum_{m=1}^{k-1}(nm+1)(-n)^{k-m-1} A_m A_{k-m}\eqno (1.3)
$$
\emph{Proof}: Take the recurrence relation for the coefficients of $f$, obtained in the previous proposition, and exclude the terms, whose indices are not divisible by $n$. Substituting $p^n=-n$ and dividing by $-n$, we obtain the desired result. \qed\\
~\\
So the first few terms in the expansion look like this:
\begin{align*}
&\mathfrak{T}f(x)=\frac{f(px)}{p},~p^n=-n ~\Rightarrow \\
&\Rightarrow f(x)=x+Ax^{n+1}+\frac{n+1}{n+2}A^2x^{2n+1}+\frac{(n+1)(n^2-n-1)}{(n+2)(n^2-3)}A^3 x^{3n+1}+\\
&+\frac{(n+1)(n^6-2n^4+4n^3-n^2-6n-2)}{(n+2)^2 (n^2-3)(n^3+4)}A^4 x^{4n+1}+...
\end{align*}
for some $A$. Let us denote this expression by $f_{n}(A^{\nicefrac{1}{n}}x)/A^{\nicefrac{1}{n}}$.\\

Since in case $n=1$, this series turns into $(e^{2Ax}-1)/2A$, it is interesting to look at coefficients of $x/f(x)$, as a sort of generalization of the generating function of Bernoulli numbers.\\
~\\
\textbf{Proposition 1.4}. Suppose $A^{\nicefrac{1}{n}}x/f_{n}(A^{\nicefrac{1}{n}}x)=\sum_{k=0}^{\infty} B_k x^{nk}$. Then $B_1=-A$, and the following relation holds true:\\
$$
(k-(-n)^{k-1})B_k = \sum_{m=1}^{k-1}(-n)^{m-1} B_m B_{k-m}\eqno (1.4)
$$
\emph{Proof}: The equality for $B_1$ is obvious. Taking the action of $-x^2\frac{\mathrm{d}}{\mathrm{d}x}\frac{1}{x}$ on the series $x/f(x)$, we obtain the recurrence relation.\qed\\
~\\
So the first few terms look like:
\begin{align*}
\frac{A^{\nicefrac{1}{n}}x}{f_{n}(A^{\nicefrac{1}{n}}x)}=&1-Ax^n+\frac{A^2 x^{2n}}{n+2}-\frac{A^3(n-1)x^{3n}}{(n+2)(n^2-3)}+\frac{A^4(n^4-n^2+4n-2)x^{4n}}{(n+2)^2(n^2-3)(n^3+4)}-\\
&-\frac{A^5(n-1)(n^2+2)(n^4-n^2+3n-1)x^{5n}}{(n+2)^2(n^2-3)(n^3+4)(n^4-5)}+...
\end{align*}
One may notice, that altought it is quite easy to calculate the coefficients of considered series, it is hard to recognise any pattern in their behaviour (compared to the case $n=1$). In order to promote greater understanding of the operator $\mathfrak{T}$, we suggest to consider the action of its conjugate $\mathrm{Q}\mathfrak{T}\mathrm{Q}$, where by $\mathrm{Q}$ we denote the operator $(\cdot)^{inv}$ $(f(f^{inv}(x))=x)$ for conciseness of notation. It should be noted, that $-\ln(1-x)$ is the image of \mbox{$\ln(1+x)$} under its action, and $\arctan(x)$ is the image of $\arcsin(x)$. Such an approach allows us to link together the solutions of equations $\mathfrak{T}f(x)=\frac{1}{p}f(px)$ and $\mathfrak{T}\mathrm{Q}f(x)=\frac{1}{p}f(px)$, \mbox{and here is why.}
\newpage
\noindent\textbf{Observation 1.1}:
$$
\begin{CD}
\displaystyle\int_{0}^{x}\frac{f^{inv}(t)}{t}dt @>\mathrm{Q}\mathfrak{T}\mathrm{Q}>>\displaystyle\int_{0}^{x}\frac{t}{\mathfrak{T}f(t)}dt
\end{CD}\eqno (1.5)
$$
\emph{Proof}: Suppose $\mathfrak{T}g_1=g_1/(g_1)'=g_2$ and $\varphi_1=(g_1)^{inv}$, $\varphi_2=(g_2)^{inv}$, i.e.
$$
\begin{CD}
\displaystyle\varphi_1 @>\mathrm{Q}\mathfrak{T}\mathrm{Q}>>\displaystyle\varphi_2
\end{CD}
$$
Then it follows that $x(\varphi_1)'=g_2(\varphi_1) \Rightarrow (x(\varphi_1)')^{inv}=g_1(\varphi_2)$. Hence $\mathfrak{T} (x(\varphi_1)')^{inv}=$ $=g_2(\varphi_2)/(\varphi_2)'=x/(\varphi_2)' \Rightarrow (\varphi_2)'=x/\mathfrak{T} (x(\varphi_1)')^{inv}$\qed\\
~\\
\textbf{Observation 1.2}: If $p_n(\alpha)$ is a sequence of polynomials of binomial type, then for any $\alpha$:
$$
\begin{CD}
\displaystyle\int_{0}^{x}\sum_{n=0}^{\infty} \frac{p_n(\alpha)t^n}{n!}dt @>\mathrm{Q}\mathfrak{T}\mathrm{Q}>>\displaystyle\int_{0}^{x}\sum_{n=0}^{\infty} \frac{p_n(-n\alpha)t^n}{n!}dt 
\end{CD}\eqno (1.6)
$$
This statement immediately follows from the Lagrange inversion theorem \cite{Lag}:
$$
(\mathfrak{T}f)^{inv}=\sum_{n=1}^{\infty} \frac{x^n}{n}\res_{t=0}\frac{(f')^n}{f^n}dt=\sum_{n=1}^{\infty} \frac{x^n}{n}\res_{t=0}\frac{[(f^{inv})']^{1-n}}{t^n}dt
$$
And so this statement is true not only for polynomials of binomial type, generated by $e^{\alpha \varphi(x)}$, but for any sequence of polynomials, generated by $(\varphi(x)/x)^{\alpha}$ (here $\varphi(x) \in x+x^2\mathbb{C}[[x]]$).\\
~\\
\textbf{Remark 1.1}. In the previous paper we considered the action of operator $\mathrm{Q}e^{-x}\mathrm{Q}$, taking into account how the polynomials of binomial type change under its action on corresponding functions (we recall that the polynomials transform like $p_n(\alpha) \rightarrow \frac{\alpha}{\alpha+n} p_n(\alpha+n)$, and so in particular case of the image of $f^{inv}(x)=x$ we end up with Abel polynomials \cite{Abl}). Then, according to the latter observations, the following relations are equivalent:

$$
\begin{CD}
    \varphi(x) @>\mathrm{Q}e^{-x}\mathrm{Q}>>\gamma(x)\\
    @.                                                                           @.\\
    xe^{-\varphi(x)} @>\mathrm{Q}>>xe^{\gamma(x)}\\
    @.                                                                           @.\\
    \displaystyle\int_{0}^{x}e^{-\varphi(t)}dt @>\mathrm{Q}\mathfrak{T}\mathrm{Q}>>\displaystyle\int_{0}^{x}1+t\gamma'(t)dt \\
\end{CD}
$$
~\\
~\\
\textbf{Observation 1.3} (\emph{conjugation} $x^{\nicefrac{1}{n}}\circ (\cdot)\circ x^n$): Note that for the formal power series $f(x) \in x+x^2 \mathbb{C}[[x]]$ the conjugation $x^{\nicefrac{1}{n}}\circ f(x)\circ x^n$: $f(x) \rightarrow x (f(x^n)/x^n)^{\nicefrac{1}{n}}$ is well-defined. For the simplicity we denote the result by $f^{\nicefrac{1}{n}}(x^n)$. Then the following holds:
\begin{align*}
\mathfrak{T}&(f^{\nicefrac{1}{n}}(x^n))=(\mathfrak{T}f)(x^n)\cdot x^{1-n}\\
\mathrm{Q}&(f^{\nicefrac{1}{n}}(x^n))=(\mathrm{Q}f)^{\nicefrac{1}{n}}(x^n) \text{, and hence}\\
\tag{1.7}\mathfrak{T}\mathrm{Q}&(f^{\nicefrac{1}{n}}(x^n))=(\mathfrak{T}\mathrm{Q}f)(x^n)\cdot x^{1-n}
\end{align*}
Such an observation allows us to rewrite the property $\mathfrak{T}f(x)=\frac{1}{p}f(px)$ in terms of inverse functions. Suppose $p^n=-n$ and the first coefficient of $f$ is equal to $A$. Then the series of $f^{inv}$ has the following form:
$$
f^{inv}(x)=x-Ax^{n+1}+\frac{(A(n+1))^2 x^{2n+1}}{n+2}+\frac{(A(n+1))^3(2+4n-3n^2)x^{3n+1}}{2(n+2)(n^2-3)}+...
$$
which in case $n=1$ turns into $\ln(1+2Ax)/2A$. For convenience we put $A^{-1}=n(n+1)$ and rewrite it in the form:
$$
-n\ln((f^{inv})')(x^{\nicefrac{1}{n}})=x-\frac{3x^2}{2(n+2)}+\frac{x^3(17n-29)}{6(n+2)(n^2-3)}-... \in  x+x^2 \mathbb{C}[[x]]
$$
Let us denote the resulting series by $\varphi_n(x)$. Then for general constant $A$ the \mbox{following holds:}
$$
(f_n(A^{\nicefrac{1}{n}}x)/A^{\nicefrac{1}{n}})^{inv}=\int_{0}^{x} e\strut^{-\tfrac{1}{n}\varphi_n(An(n+1)t^n)}dt \eqno(1.8)
$$
Thus the property $\mathrm{Q}\mathfrak{T}\mathrm{Q}~f^{inv}(x)=\frac{1}{p}f^{inv}(px)$ implies that
$$
\begin{CD}
\displaystyle\int_{0}^{x}e\strut^{-\tfrac{1}{n}\varphi_n(t^n)}dt @>\mathrm{Q}\mathfrak{T}\mathrm{Q}>>\displaystyle\int_{0}^{x}e\strut^{-\tfrac{1}{n}\varphi_n(-nt^n)}dt 
\end{CD}
$$
And hence, according to (1.5), this property can be rewritten in the following form:
$$
\mathfrak{T}\mathrm{Q}~xe\strut^{-\tfrac{1}{n}\varphi_n(x^n)}=xe\strut^{\tfrac{1}{n}\varphi_n(-nx^n)}
$$
or, taking (1.7) into account:
$$
\boxed{\mathfrak{T}\mathrm{Q}~xe\strut^{-\varphi_n(x)}=xe\strut^{\tfrac{1}{n}\varphi_n(-nx)}}\eqno(1.9)
$$
Note that this relation can be considered as a definition of the series $\varphi_n(x)\in x+x^2\mathbb{C}[[x]]$.\\
~\\
\textbf{Remark 1.2}. It should be noticed that the existense of such a series $\varphi_n(x)$ does not necessarily mean, that $n$ must be a natural number. The condition $n \in \mathbb{N}$ is used only when one wants, going in the opposite direction, to construct the meaningful series $f_n(x)$, which consists of the powers $x^{dn+1}$. The equation for $f$ can be rewritten in another form, according to Observation 1.3, as follows: $\mathfrak{T}~f^n(x^{\nicefrac{1}{n}})=xf(px^{\nicefrac{1}{n}})/(px^{\nicefrac{1}{n}})$. That means with the solution of equation $\mathfrak{T}~x\exp(Mq_1(x))=x\exp(Mq_1(-px)/p)$ at our disposal, we can construct the solution of equation $\mathfrak{T}\mathrm{Q}~x\exp(-Mq_2(x))=x\exp(Mq_2(-px)/p)$ (in this case, however, $q_2$ may not be the invertible series, being of the form $x^k+x^{k+1}\mathbb{C}[[x]]$). Let us do it without using the trick with conjugation by the powers of $x$.\\
~\\
\textbf{Proposition 1.5}. Suppose for some $\sigma, q \in x+x^2\mathbb{C}[[x]], M,p \in \mathbb{C}$ we have
$$
\begin{CD}
    xe\strut^{Mq(x)} @>\mathfrak{T}>>xe\strut^{\tfrac{M}{p}q(-px)} , ~~~~~@. xe\strut^{Mq(x)} @>\mathrm{Q}>>xe\strut^{-M\sigma(x)} \\
\end{CD}
$$
Then the following holds true:
$$
\begin{CD}
    x(1-Mx\sigma'(x))\strut^p e\strut^{-M\sigma(x)} @>\mathfrak{T}\mathrm{Q}>>x(1+Mpx\sigma'(-px))\strut^{-1} e\strut^{\tfrac{M}{p}\sigma(-px)}\\
\end{CD}
$$
\emph{Proof}: Obviously,
$$ [(x\exp(-M\sigma(x)))']^p=(1-Mx\sigma'(x))^p \exp(-Mp\sigma(x)) $$
In addition, from the assumption it follows that:
\begin{align*}
&\sigma(x \exp(Mq(x)))=q(x)\\
&q(x \exp(-M\sigma(x)))=\sigma(x)~~\text{and}\\
&x(1+Mxq'(x))^{-p}=x\exp(Mq(-px))
\end{align*}
\newpage
\noindent Hence
\begin{align*}
&\mathfrak{T}\mathrm{Q}\left(x\exp(-M\sigma(x))(1-Mx\sigma'(x))\strut^p\right)=\\
&=\mathfrak{T}\mathrm{Q}\left(x\exp((p-1)M\sigma(x))[(x \exp(-M\sigma(x)))']\strut^p\right)=\\
&=\mathfrak{T}\mathrm{Q}\left[\left(x\exp((p-1)M\sigma(x))[(x \exp(-M\sigma(x)))']\strut^p\right)\circ xe\strut^{Mq(x)} \circ xe\strut^{-M\sigma(x)}\right]=\\
&=\mathfrak{T}\mathrm{Q}\left[\left(x\exp(Mpq(x))[(x \exp(Mq(x)))']\strut^{-p}\right)\circ xe\strut^{-M\sigma(x)}\right]=\\
&=\mathfrak{T}\mathrm{Q}\left[x\exp(Mq(-px))\circ xe\strut^{-M\sigma(x)}\right]=\mathfrak{T}\left[xe\strut^{Mq(x)}\circ(x\exp(Mq(-px)))^{inv}\right]=\\
&=\mathfrak{T}\left[xe\strut^{Mq(x)}\circ xe\strut^{-M\sigma(-px)}\right]=\frac{1}{\left(xe\strut^{-M\sigma(-px)}\right)'}\left[\left(\mathfrak{T}xe\strut^{Mq(x)}\right)\circ xe\strut^{-M\sigma(-px)}\right]=\\
&=\frac{1}{\left(xe\strut^{-M\sigma(-px)}\right)'}\left[xe\strut^{\tfrac{M}{p}q(-px)}\circ xe\strut^{-M\sigma(-px)}\right]=\frac{x\exp\left(\tfrac{M(1-p)}{p}\sigma(-px)\right)}{\left(x\exp(-M\sigma(-px))\right)'}=\\
\tag*{\qed}&=x(1+Mpx\sigma'(-px))\strut^{-1} \exp\left(\frac{M}{p}\sigma(-px)\right)
\end{align*}
That means the equation, which we are interested in, originates from more general one for $p=n$:
$$
\begin{CD}
    xe\strut^{Mq(x)} @>\mathfrak{T}>>xe\strut^{\tfrac{M}{n}q(-nx)} ~\Rightarrow~ @. xe\strut^{\tfrac{M}{n}q(x^n)} @>\mathfrak{T}>>xe\strut^{\tfrac{M}{n}q(-nx^n)} \\
\end{CD}
$$
Thus the existence of companion equation (1.9) for the series $\varphi_n$ is the direct corollary of this proposition. In particular, that means we can consider the general series $\varphi_p \in x+x^2\mathbb{C}[[x]]$, whenever the equation
$$
\mathfrak{T}\mathrm{Q}~xe\strut^{-\varphi_p(x)}=xe\strut^{\tfrac{1}{p}\varphi_p(-px)}
$$
has the solution for given $p$. At this point the only case $\varphi_1(x)=\ln(1+x)$ is known. Let us consider some other interesting cases.

\begin{center}
	\textbf{The equation} $\mathfrak{T}^2f(x)=\dfrac{f(2x)}{2}$
\end{center}
As it was mentioned in the previous paper, the equation $\mathfrak{T}^2f(x)=f(2x)/2$ has two elementary solutions: $\sin x$ and $\tan x$. In this paragraph we provide the general solution of this equation.\\
~\\
\underline{\textbf{Theorem 1.1}}. Let $\Theta_k(x)$ be the Jacobi elliptic function ($\mathrm{sn}(k | x)$ in conventional notation), i.e.
$$
\Theta_k^{inv}(x)=\int_{0}^{x} \frac{dt}{\sqrt{1-t^2}\sqrt{1-k^2 t^2}}
$$
Then it follows that
$$
\mathfrak{T}^2 \Theta_k(x)=\frac{\Theta_k(2x)}{2} \eqno(1.10)
$$
\emph{Proof}: By definition,
$$
\Theta_k'(x)^2 = (1-\Theta_k^2(x))(1-k^2\Theta_k^2(x)) ~\Rightarrow~ \Theta_k''(x)=2k^2\Theta_k^3(x)-(k^2+1)\Theta_k(x)
$$
Also, the Jacobi elliptic function obeys the addition theorem (see \cite{Jac}):
$$
\Theta_k(x+y)=\frac{\Theta_k(x)\Theta_k'(y)+\Theta_k'(x)\Theta_k(y)}{1-k^2\Theta_k^2(x)\Theta_k^2(y)}
$$
Hence
\begin{align*}
\mathfrak{T}^2\Theta_k(x)=&\frac{\Theta_k(x)\Theta_k'(x)}{\Theta_k'(x)^2-\Theta_k(x)\Theta_k''(x)}=\\
=&\frac{\Theta_k(x)\Theta_k'(x)}{1-(k^2+1)\Theta_k^2(x)+k^2\Theta_k^4(x)-2k^2\Theta_k^4(x)+(k^2+1)\Theta_k^2(x)}=\\
\tag*{$\blacksquare$}=&\frac{\Theta_k(x)\Theta_k'(x)}{1-k^2\Theta_k^4(x)}=\frac{\Theta_k(2x)}{2}
\end{align*}
\textbf{Corollary 1.1}: $\forall A \in \mathbb{C}$ the series $\Theta_k(Ax)/A$ is also a solution of this equation.\\
~\\
Notice that we don't really need to know the whole addition formula to prove this theorem, the duplication formula only is required. Thus we can use another approach to prove this theorem. Consider the chain, whose main properties were overviewed in the previous paper (we use the notation $\Delta_p=(e^{px}-1)/p$, $y_p=\Delta_p e^{-x}$):
$$
\begin{CD}
    y_p @>\mathfrak{T}>> \displaystyle\frac{\Delta_{-p}}{1-\Delta_{-p}} @>\mathfrak{T}>> \Delta_{p}(1-\Delta_{-p}) @>\mathfrak{T}>> \displaystyle\frac{\mathfrak{T}y_p}{1+(p-1)(\mathfrak{T}y_p)^2}\\
\end{CD}
$$
Now notice the transformation of inverse series, corresponding to the links $\mathfrak{T}(\mathfrak{T}y_p)=\mathfrak{T}^2y_p$:
$$
\begin{CD}
\displaystyle\int_{0}^{x}\frac{dt}{1+(2-p)t+(1-p)t^2} @>\mathrm{Q}\mathfrak{T}\mathrm{Q}>>\displaystyle\int_{0}^{x}\frac{dt}{\sqrt{1+2(p-2)t+p^2t^2}}
\end{CD}
$$
Then after conjugation $t^{\nicefrac{1}{2}}\circ(\cdot)\circ t^2$ one can obtain the transformation:
$$
\begin{CD}
\displaystyle\int_{0}^{x}\frac{dt}{\sqrt{1+(2-p)t^2+(1-p)t^4}} @>\mathrm{Q}\mathfrak{T}\mathrm{Q}>>\displaystyle\int_{0}^{x}\frac{dt}{\sqrt{1+2(p-2)t^2+p^2t^4}}
\end{CD}
$$
which, obviously, goes on in the following way:
$$
\begin{CD}
@>\mathrm{Q}\mathfrak{T}\mathrm{Q}>>\displaystyle\int_{0}^{x}\frac{dt}{\sqrt{1+4(2-p)t^2+16(1-p)t^4}}
\end{CD}
$$
That means the only thing required to prove the theorem is to notice, how the transformation $\mathfrak{T}(\Delta_{-p}/(1-\Delta_{-p}))=\Delta_{p}(1-\Delta_{-p})$ looks like in terms of inverse series.\\
~\\
\textbf{Remark 1.3}. It should also be mentioned that the addition formula for the redefined function
$$
\Theta_p^{inv}(x)=\int_{0}^{x} \frac{dt}{\sqrt{1+(2-p)t^2+(1-p) t^4}}~~ (1-p=k^2, t\rightarrow it)
$$
can be rewritten in the following way:
$$
\Theta_p(x+y)=\left(\frac{1}{\mathfrak{T}\Theta_p}(x)+\frac{1}{\mathfrak{T}\Theta_p}(y)\right)\cdot\left[(\mathfrak{T}^3y_p)\circ (\mathfrak{T}y_p)^{inv}\right](\Theta_p(x)\Theta_p(y))
$$
~\\
\underline{\textbf{Theorem 1.2}}. (\emph{uniqueness})
$$
\mathfrak{T}^2f(x)=\frac{f(2x)}{2} ~\Rightarrow~ \exists A,k \in \mathbb{C}: f(x)=\frac{\Theta_k(Ax)}{A} \eqno(1.11)
$$
The proof is divided in a few parts.
Let us denote by $\left.f\right|_n$ the polynomial $x+f_1 x^2+...+$ $+f_{n-1}x^n$, corresponding to the series $f(x)$, reduced to the $n$-th term.

\newpage
\noindent\textbf{Observation 1.4}: $\left.(f+qx^{n+1})\right|_{n+1}=\left.(\left.f\right|_{n+1}e^{qx^n})\right|_{n+1}$\\
~\\
\textbf{Observation 1.5}: $\left.(\mathfrak{T}(\left.f\right|_n))\right|_n=\left.(\mathfrak{T}f)\right|_n$\\
~\\
\textbf{Proposition 1.6}.
$$
(\mathfrak{T}\left.(\left.f\right|_{n+1}e^{qx^n}))\right|_{n+1}=\left.(\left.(\mathfrak{T}f)\right|_{n+1}e^{-nqx^n})\right|_{n+1}\eqno(1.12)
$$
\emph{Proof}:
\begin{align*}
(\mathfrak{T}\left.(\left.f\right|_{n+1}e^{qx^n}))\right|_{n+1}=&\cfrac{1}{\cfrac{1}{\mathfrak{T}(\left.f\right|_{n+1})}+nqx^{n-1}}\left.\strut\right|_{n+1}=\\
=&[\mathfrak{T}(\left.f\right|_{n+1})](1+nqx^{n-1}\mathfrak{T}(\left.f\right|_{n+1}))^{-1}\left.\strut\right|_{n+1}=\\
=&[\mathfrak{T}(\left.f\right|_{n+1})](1-nqx^n)\left.\strut\right|_{n+1}=\\
\tag*{\qed}=&\left.(\left.(\mathfrak{T}f)\right|_{n+1}e^{-nqx^n})\right|_{n+1}
\end{align*}
(note: this statement is nothing more than the brief reformulation of the following two facts: if $f(x)\approx x+f_1x^2+...+f_k x^{k+1}$, then the polynomials $A_m(f_1,...,f_m)$ in expansion $\mathfrak{T}f(x)\approx x+A_1x^2+...+A_k x^{k+1}$ indeed do not depend on the higher-order coefficients $f_{m+1},f_{m+2},...$, and that the polynomial $A_k$ depends on $f_k$ linearly with coefficient $-k$, i.e. $A_k(f_1,...,f_k+q)=A_k(f_1,...,f_k)-kq$)\\
~\\
\textbf{Proposition 1.7}. Suppose $\mathfrak{T}^2 f(x)=f(2x)/2$. Then $f(x)=-f(-x)$.\\
\emph{Proof}: Suppose such $f(x)$ is nontrivial, i.e. there exists the least $n \geqslant 1$, such that the coefficient $\alpha$ of $x^{n+1}$ is not equal to zero. Then
$$
\left.f\strut\right|_{n+1}=x+\alpha x^{n+1} ~\Rightarrow~ \left.(\mathfrak{T}^2f)\right|_{n+1}=x+n^2\alpha x^{n+1}
$$
On the other hand, we have $(\mathfrak{T}^2f)|_{n+1}=x+2^n\alpha x^{n+1}$. Hence $2^n=n^2 ~\Rightarrow~ n=2$ or $n=4$. Consider only the first case without loss of generality, i.e. $\left.f\strut\right|_{3}=x+\alpha x^3$, where now $\alpha$ can be equal to zero. Now if there exist some $A_i$ and some $q$, such that
$$
\left.f\strut\right|_{2n+2}=x+A_1 x^3+A_2 x^5+...+A_n x^{2n+1}+qx^{2n+2}
$$
Then
$$
\left.(\mathfrak{T}^2f)\strut\right|_{2n+2}=x+2^2A_1 x^3+2^4A_2 x^5+...+2^nA_n x^{2n+1}+(2n+1)^2qx^{2n+2}
$$
according to Proposition 1.6 and the fact that the coefficients $A_i$ of odd powers $x^{2n+1}$ can not give contributions to the coefficient of even power $x^{2n+2}$ under the action of $\mathfrak{T}$. But it means that $(2n+1)^2q=2^{2n+1}q$, and since the equation $(2n+1)^2=2^{2n+1}$ does not have any integer solutions, we have $q=0$. It means that all of the coefficients of even powers in the expansion of $f$ are equal to $0$.\qed\\
~\\
\textbf{Proposition 1.8}. Suppose $f$ is the solution, which has the expansion $f|_{5}=x+\alpha x^3+\beta x^5$. Suppose $k, A$ is any pair of solutions of the system of equations $-A^2(k^2+1)=3!\alpha$, $A^4(1+14k^2+k^4)=5!\beta$. Then $f(x)=\Theta_k(Ax)/A$.\\
\emph{Proof}: One can easily derive the first few terms of Jacobi function from its definition and obtain the expansion
$$
\frac{\Theta_k(Ax)}{A} \approx x-\frac{k^2+1}{6}A^2 x^3+\frac{k^4+14k^2+1}{5!}A^4 x^5
$$
That means it coincides with the expansion of $f$ right up to the fifth term. Now suppose $q$ is the first deviation between the expansions of some term $x^{2n+1}$, i.e.
$$
\left.f\strut\right|_{2n+1}=\frac{\Theta_k(Ax)}{A}\left.\strut\right|_{2n+1}+qx^{2n+1}
$$
Then, according to Proposition 1.6,
$$
\left.(\mathfrak{T}^2 f)\strut\right|_{2n+1}=\left(\mathfrak{T}^2\frac{\Theta_k(Ax)}{A}\right)\left.\strut\right|_{2n+1}+(2n)^2qx^{2n+1}
$$
So either $(2n)^2=2^{2n}$, or $q=0$. But we have $n>2$ due to the choice of \mbox{$k, A \Rightarrow (2n)^2 \neq 2^{2n}$}, hence $q=0$.\qed\\ That means the series $\Theta_k(Ax)/A$ are the only series of the form $\in x+x^2\mathbb{C}[[x]]$, satisfying the equation $\mathfrak{T}^2f=f(2x)/2$.\hfill\ensuremath{\blacksquare}
\begin{center}
	$\varphi_n$ \textbf{in case} $n \neq 1$
\end{center}
Consider again the series $\varphi_n \in x+x^2\mathbb{C}[[x]]$, defined as the solutions of equation
$$\mathfrak{T}\mathrm{Q}~xe\strut^{-\varphi_n(x)}=xe\strut^{\tfrac{1}{n}\varphi_n(-nx)}$$
Since we now know how to solve the equation $\mathfrak{T}^2 f(x)=f(2x)/2$, we are ready to write down the solutions of equation $\mathfrak{T}f(x)=f(px)/p$ in cases $p^2=-2$ and $p^4=-4$. For that purpose we note first that the following relation holds (see the comment after \mbox{Corollary 1.1):}
$$
\mathfrak{T}\Theta_k(x)=\frac{1}{i(k+1)}\Theta_{\frac{k-1}{k+1}}(i(k+1)x)
$$
Then in case $k=\sqrt{2}-1$, since $\Theta_k=\Theta_{-k}$, we have
$$\mathfrak{T}\Theta_{\sqrt{2}-1}(x)=\frac{1}{i\sqrt{2}}\Theta_{1-\sqrt{2}}(i\sqrt{2}x)=\frac{1}{\sqrt{-2}}\Theta_{\sqrt{2}-1}(\sqrt{-2}x)$$
Similarly, in case $k=i$, since the lemniscatic elliptic function
$$\Theta_i(x)=\left(\displaystyle\int_{0}^{x}\frac{dt}{\sqrt{1-t^4}}\right)^{inv}$$
is invariant under the transformation $x \rightarrow ix$, $\Theta_i(ix)=i\Theta_i(x)$, we have
$$\mathfrak{T}\Theta_{i}(x)=\frac{1}{i(1+i)}\Theta_{i}(i(1+i)x)=\frac{1}{1+i}\Theta_{i}((1+i)x)=\frac{1}{\sqrt[4]{-4}}\Theta_{i}(\sqrt[4]{-4}x)$$
Since in terms of $f_n$ we have the following expansions:
\begin{align*}
&\frac{1}{A^{\nicefrac{1}{2}}}\left.f_2(A^{\nicefrac{1}{2}}x)\strut\right|_{3}=x+Ax^3, &&\frac{1}{p}\Theta_{\sqrt{2}-1}(px)\left.\strut\right|_{3}=x-p^2\frac{(\sqrt{2}-1)^2+1}{6}x^3\\
&\frac{1}{A^{\nicefrac{1}{4}}}\left.f_4(A^{\nicefrac{1}{4}}x)\strut\right|_{5}=x+Ax^5, &&\frac{1}{p}\Theta_i(px)\left.\strut\right|_{5}=x-\frac{p^4}{10}x^5
\end{align*}
after substitution $-p_1^2((\sqrt{2}-1)^2+1)=6$ and $-p_2^4=10$, we can now show that
\begin{align*}
&f_2(x)=\frac{1}{p_1}\Theta_{\sqrt{2}-1}(p_1 x), && f_4(x)=\frac{1}{p_2}\Theta_i(p_2 x)
\end{align*}
Or, in other words, we have the following expressions:
\begin{align*}
&f_2(x)^{inv}=\int_{0}^{x} \frac{dt}{\sqrt{1+6t^2+\frac{9}{2}t^4}}; && f_4(x)^{inv}=\int_{0}^{x}\frac{dt}{\sqrt{1+10t^4}}
\end{align*}
from which, according to (1.8), the identities
\begin{align*}
&\boxed{\varphi_2(x)=\ln\left(1+x+\frac{x^2}{8}\right)} &&\boxed{\varphi_4(x)=2\ln\left(1+\frac{x}{2}\vphantom{\frac{x^2}{8}}\right)} 
\end{align*}
follow. And indeed, we have:
\begin{align*}
&\mathfrak{T}\mathrm{Q}~xe\strut^{-\varphi_2(x)}=\mathfrak{T}\left(\frac{x}{1+x+\tfrac{1}{8}x^2}\right)^{inv}=\mathfrak{T}~\frac{2x}{1-x+\sqrt{1-2x+\frac{1}{2}x^2}}=x\sqrt{1-2x+\frac{1}{2}x^2}=\\
&=xe\strut^{\tfrac{1}{2}\varphi_2(-2x)}~;\\
~\\
&\mathfrak{T}\mathrm{Q}~xe\strut^{-\varphi_4(x)}=\mathfrak{T}\left(\frac{x}{\left(1+\tfrac{1}{2}x\right)^2}\right)^{inv}=\mathfrak{T}~\frac{2x}{1-x+\sqrt{1-2x}}=x\sqrt{1-2x}=xe\strut^{\tfrac{1}{4}\varphi_4(-4x)}
\end{align*}
We note that besides the considered natural cases, there are fractional ones, which remain to be explored. And although there is no apparent structure in these cases, there are actually two additional values of $n$, for which the series $\varphi_n$ can be written down explicitly. We will return to this topic soon. For now consider the following two observations.\\
\textbf{Observation 1.6}: Suppose, the equation $\mathfrak{T}\mathrm{Q}~xe\strut^{-\varphi_\nu(x)}=xe\strut^{\tfrac{1}{\nu}\varphi_\nu(-\nu x)}$ is solvable for some $\nu\in\mathbb{C}$, and the series $\varphi_\nu \in x+x^2\mathbb{C}[[x]]$ is the corresponding solution. Then if $\varphi_\nu$ generates polynomials of binomial type
$$
e^{\alpha\varphi_\nu(x)}=\sum_{n=0}^\infty \frac{p_n(\alpha)x^n}{n!}
$$
the following holds true:
$$
(-\nu)^n p_n\left( -\frac{1}{\nu}\right)=p_n(n)\eqno(1.13)
$$
(see Observation 1.2)\\
\textbf{Observation 1.7}: The equation $\mathfrak{T}\mathrm{Q}~xe\strut^{-\varphi_\nu(x)}=xe\strut^{\tfrac{1}{\nu}\varphi_\nu(-\nu x)}$ is equivalent to the equation
$$
e\strut^{-\tfrac{1}{\nu}\varphi_\nu(-\nu x)}=\sum_{n=0}^{\infty} \frac{\left(x\tfrac{d}{dx}\right)^n}{n!}(\varphi_\nu(x))^n\eqno(1.14)
$$
\emph{Proof}: This statement is the direct corollary of the latter observation. Taking into account the equality $p_n(\alpha)=n!\res_{t=0}e^{\alpha\varphi_\nu(t)}t^{-n-1}dt$, we obtain
\begin{align*}
\sum_{n=0}^\infty &\frac{p_n(n) x^n}{n!}=\sum_{n=0}^\infty \frac{x^n}{n!}\sum_{k=0}^n \frac{n^k n!}{k!}\res_{0}\frac{\varphi_\nu(t)^k}{t^{n+1}}dt=\sum_{k=0}^\infty\sum_{n=k}^\infty n^kx^n\res_{0}\frac{\varphi_\nu(t)^k}{t^{n+1}}dt\\
\tag*{\qed}&\Rightarrow\sum_{n=0}^\infty \frac{(-\nu)^n p_n\left( -\tfrac{1}{\nu}\right)x^n}{n!}=\sum_{n=0}^\infty \frac{p_n(n) x^n}{n!}=\sum_{n=0}^{\infty} \frac{\left(x\tfrac{d}{dx}\right)^n}{n!}(\varphi_\nu(x))^n
\end{align*}
\textbf{Proposition 1.9}. Consider the series $\psi(x)=\left(x e^{\mathrm{Ei}(x)-\ln|x|-\gamma}\right)^{inv}$, introduced in the previous paper, which is equivalently determined by the equation $x\psi'(x)=\psi(x) e^{-\psi(x)}$. Consider, in addition, Lambert $W$-function $W(x)=(xe^x)^{inv}$. Then the following limits exist (in the sense that for each $k$, the sequence of coefficients of $x^k$ converges):
\begin{align*}
\lim_{p\to 0}~&\varphi_p(x)=-\ln\frac{\psi(x)}{x} &&\lim_{p\to \infty}\frac{1}{p}\varphi_p(px)=\ln(1+W(x))\\
\lim_{p\to \infty}&\varphi_p(x)=x &&\lim_{p\to 0}~\frac{1}{p}\varphi_p(px)=x
\end{align*}
And so the following identities are true:
\begin{align*}
&\boxed{\varphi_\infty(x)=x\vphantom{-\ln\frac{\psi(x)}{x}}} &&\boxed{\varphi_0(x)=-\ln\frac{\psi(x)}{x}}
\end{align*}
\emph{Proof}: According to Observation 1.7, if $\varphi_p(x)$ has the expansion $\varphi_p(x)=A_1(p)x+ +A_2(p)x^2+...$, and $A_1(p)=1$, one can write down the following relation, comparing the coefficients of LHS and RHS of (1.14):
$$
\sum_{k=1}^{n} \frac{(-p)^{n-k}}{k!}\sum_{\substack{\sum_1^k m_i=n\\ m_i>0}} A_{m_1}(p)...A_{m_k}(p)=\sum_{k=1}^{n} \frac{n^k}{k!}\sum_{\substack{\sum_1^k m_i=n\\ m_i>0}} A_{m_1}(p)...A_{m_k}(p)
$$
Or, equivalently
$$
(n-(-p)^{n-1})A_n(p)=\sum_{k=2}^{n} \frac{(-p)^{n-k}-n^k}{k!}\sum_{\substack{m_1+...+m_k=n\\ m_i>0}} A_{m_1}(p)...A_{m_k}(p)
$$
That means, if $p\neq -2, p^2 \neq 3, p^{3}\neq -4,...$, then the condition $A_1(p)=1$ determines the sequence $A_n(p)$. So, if $A_m(p)$ in case $p \to 0$ converge to some value $A_m(0)$ for each $m<n$, then $A_n(p)$ converges too and
$$
A_n(0)=\frac{1}{n}\sum_{k=2}^{n} \frac{\delta_{n,k}-n^k}{k!}\sum_{\substack{m_1+...+m_k=n\\ m_i>0}} A_{m_1}(0)...A_{m_k}(0)
$$
(here, as usual $\delta_{n,k}=1$, if $n=k$, and $0$ otherwise). Thus all $A_n(0)$ are well-defined. Similarly, we can rewrite this relation in the following form:
$$
(n(-p)^{1-n}-1)A_n(p)=\sum_{k=2}^{n} \frac{(-p)^{1-k}-n^k(-p)^{1-n}}{k!}\sum_{\substack{\sum_1^k m_i=n\\ m_i>0}} A_{m_1}(p)...A_{m_k}(p)
$$
Then, since $A_1(\infty)$ is defined, we can conclude that $A_n(\infty)$ are well-defined too and are in fact trivial: $\forall n>1 ~A_n(\infty)=0$. Hence,
$\lim_{p \to 0}\varphi_p(x)=\varphi_0(x)$ and $\lim_{p \to \infty}\varphi_p(x)=x$ in the sense described above. Then it follows that the two other limits $\lim_{p \to 0}\varphi_p(px)/p=x$ and $\lim_{p \to \infty}\varphi_p(px)/p$  exist, since the operator $\mathfrak{T}\mathrm{Q}$ is invertible. In case $p=0$ we end up with the equation
$$
\mathfrak{T}\mathrm{Q}~ xe\strut^{-\varphi_0(x)}=xe^{-x} ~\Rightarrow~ xe\strut^{-\varphi_0(x)}=\mathrm{Q}\mathfrak{T}^{-1}xe^{-x}~\Rightarrow~\varphi_0(x)=-\ln\frac{\psi(x)}{x}
$$
And in case $p=\infty$ we have
$$
\mathfrak{T}\mathrm{Q}~xe^{-x}=x\exp\left(\lim_{p\to \infty}\frac{1}{p}\varphi_p(-px)\right)
$$
On the other hand, the following holds true:
$$
\mathfrak{T}\mathrm{Q}~xe^{-x}=\mathfrak{T}(-W(-x))=-\frac{W(-x)}{W'(-x)}=x(1+W(-x))
$$
Hence
$$
\lim_{p\to \infty}\frac{1}{p}\varphi_p(-px)=\ln(1+W(-x)) ~\Rightarrow~ \lim_{p\to \infty}\frac{1}{p}\varphi_p(px)=\ln(1+W(x))
$$\qed\\
~\\
\textbf{Remark 1.4}. Then it follows from this proposition that the behaviour of free terms and leading coefficients of numerators and denominators of rational functions $A_n(p)$ is relatively well understood.

\begin{center}
\textbf{General comment on section 1}
\end{center}
Taking into account Remark 1.2, it is reasonable to consider the equation \mbox{$\mathfrak{T}\mathrm{Q}f(x)=\frac{1}{p}f(px)$}. It has the formal solutions whenever $p^n=n$ for some natural number $n$, to be more precise
$$\mathfrak{T}\mathrm{Q}~xe\strut^{-\tfrac{1}{n}\varphi_{-n}(x^n)}=xe\strut^{-\tfrac{1}{n}\varphi_{-n}(nx^n)}$$
for all $n \neq 2$, by the definition of $\varphi_{-n}$, according to (1.7). In case $p^2=2$ the solution takes the form $x+A_1x^5+A_2x^9+...$, i.e. it coincides with solution of the same equation in case $p^4=4$. It is interesting to know, if these series converge in some domain or they are purely formal power series.\\
Since the chain $y_p \to \mathfrak{T}y_p \to \Delta_p (1-\Delta_{-p})$ provides the major part of known examples of $\varphi_n$ for different $n$, one can consider other chains of transformations $\mathrm{Q}\mathfrak{T}\mathrm{Q}$ for different expressions, containing the functions $\gamma_p=y_p^{inv}$. An interested reader might find the main examples in Appendix A.

\newpage

\hrule

\section{Some properties of canonical continuations of polynomials of binomial type}

\hrule

~\\
~\\
In this section we continue to study the formal series $p_s(\alpha)$, which is the natural generalization of polynomials of binomial type $p_n(\alpha)$. The main topic discussed here is about the derivative $\partial/\partial s$: its connection with alternative definition of logarithm, as an ``integral'', related to the ``derivation'' $f(D)$, and how it allows to prove reflection formulae for $p_s(\alpha)$ in specific cases.\\
Here we use the following notation: the falling factorials $x(x-1)...(x-n+1)$ are denoted by $(x)_n$; an operator $\partial/\partial \alpha$ is denoted by $D$ for conciseness. Also we use the expression $f(x) \rightsquigarrow p_s(\alpha)$ or, alternatively, $p_{\vphantom{1} s}^f(\alpha)$ for $f(x) \in x+x^2\mathbb{C}[[x]]$ to denote the formal series $p_s(\alpha) \in \alpha^s+\alpha^{s-1}\mathbb{C}[[\alpha^{-1}]]$, which is determined by condition
$$
p_s(\alpha)=\alpha \left(\frac{D}{f(D)}\right)^s \alpha^{s-1}\eqno(2.1)
$$
\textbf{Remark 2.1}. This definition is consistent with expansion
$$
\sum_{n=0}^\infty \frac{p_n(\alpha)x^n}{n!}=\exp(\alpha f^{inv}(x))
$$
\textbf{Remark 2.2}. Such a notation with $D$ is just a brief reformulation of the following definition
$$
\left( \frac{x}{f(x)}\right)^s=\sum_{n=0}^\infty \frac{q_n(s)x^n}{n!} ~\Rightarrow~ p_s(\alpha)=\sum_{n=0}^\infty \binom{s-1}{n}q_n(s)\alpha^{s-n}
$$
Also the relations $f(D)p_s(\alpha)=sp_{s-1}(\alpha)$ and $\alpha\frac{1}{f'(D)}p_s(\alpha)=p_{s+1}(\alpha)$ hold.\\
\textbf{Observation 2.1}: The formal series $\frac{\partial}{\partial s}p_s(\alpha)$ is well-defined  and has expansion of the form $p_s(\alpha)\ln \alpha +r_s(\alpha)$ for some formal series $r_s(\alpha)\in \alpha^{s-1}\mathbb{C}[[\alpha^{-1}]]$ and formal element $\ln\alpha$, defined by the condition $D\alpha^q\ln\alpha=\alpha^{q-1}+q\alpha^{q-1}\ln\alpha$.\\
\emph{Proof}: The correctness of such a definition is obvious. Differentiate (2.1) by $s$ formally to obtain the following series:
\begin{align*}
\tag*{\qed}p_s(\alpha)\ln\alpha+\sum_{n=1}^\infty \alpha^{s-n}\left[\frac{\partial}{\partial s}\binom{s-1}{n}q_n(s)\right]
\end{align*}
We subsequently use the notation $\dot p_s(\alpha)$ for $\frac{\partial}{\partial s}p_s(\alpha)$.\\
~\\
\textbf{Definition}: For each $g \in x+x^2\mathbb{C}[[x]]$, define the formal element $\ln_g \alpha$ of the form $\ln\alpha+\sum_{n=1}^\infty A_n \alpha^{-n}$, which satisfies the condition $g(D)\ln_g \alpha=\alpha^{-1}$.\\
~\\
\textbf{Observation 2.2}: Suppose
$$\frac{x}{g(x)}=\sum_{n=0}^\infty \frac{A_n}{n!}x^n$$
Then it follows that
$$
\ln_g \alpha =\ln \alpha+\sum_{n=1}^\infty \frac{(-1)^{n-1}A_n}{n}\alpha^{-n}\eqno(2.2)
$$
\emph{Proof}:
\begin{align*}
\ln_g\alpha&=\frac{D}{g(D)}\ln\alpha=\sum_{n=0}^\infty \frac{A_n}{n!}D^n\ln\alpha=\ln\alpha+\sum_{n=1}^\infty \frac{A_n}{n!}D^{n-1}\alpha^{-1}=\\
\tag*{\qed}&=\ln \alpha+\sum_{n=1}^\infty \frac{(-1)^{n-1}A_n}{n}\alpha^{-n}
\end{align*}
\textbf{Proposition 2.1}. Suppose $f(x) \rightsquigarrow p_s(\alpha)$. Then
\begin{align*}
\tag{2.3}&\dot p_0(\alpha)=\ln_{\mathfrak{T}f}\alpha\\
\tag{2.4}&\dot p_1(\alpha)=\alpha \ln_f \alpha
\end{align*}
\emph{Proof}: $\left(x/f(x)\right)^s=\sum_{n\geqslant 0} q_n(s)x^n/n! \Rightarrow \forall n>0, ~q_n(0)=0$. Then it follows that
\begin{align*}
\left.\frac{\partial}{\partial s}p_s(\alpha)\right|_{s=0}&=\ln\alpha+\sum_{n=1}^\infty \binom{-1}{n}\alpha^{-n}\left.\frac{\partial}{\partial s}q_n(s)\right|_{s=0}=\\
&=\ln\alpha+\sum_{n=1}^\infty (-1)^n\alpha^{-n}\left.\frac{\partial}{\partial s}q_n(s)\right|_{s=0}=\ln\alpha+\sum_{n=1}^\infty (-\alpha)^{-n}q_n'(0)
\end{align*}
On the other hand,
$$\sum_{n=0}^\infty \frac{q_n'(0)x^n}{n!}=\ln\frac{x}{f(x)} ~\Rightarrow~ \frac{xf'(x)}{f(x)}=1-\sum_{n=1}^\infty \frac{q_n'(0)x^n}{(n-1)!}$$
Then, according to (2.2), the following holds:
$$
\ln_{\mathfrak{T}f}\alpha=\ln\alpha+\sum_{n=1}^\infty \frac{(-1)^{n-1}}{n}\alpha^{-n}(-nq_n'(0))=\dot p_0(\alpha)
$$
Similarly,
$$
\dot p_1(\alpha)=\alpha\ln\alpha+\sum_{n=1}^\infty q_n(1)\alpha^{1-n}\left.\frac{\partial}{\partial s}\binom{s-1}{n}\right|_{s=1}=\alpha\ln\alpha-\alpha\sum_{n=1}^\infty \frac{q_n(1)}{n}(-\alpha)^{-n}
$$
And that is indeed $\alpha \ln_f \alpha$.\qed\\
~\\
\textbf{Remark 2.3}. It should be mentioned that from the relation $f(D)p_s(\alpha)=sp_{s-1}(\alpha)$ it follows that $f(D)\dot p_s(\alpha)=p_{s-1}(\alpha)+s\dot p_{s-1}(\alpha)$. Substituting $s=0$, we obtain the relation $f(D)\dot p_0(\alpha)=p_{-1}(\alpha)=f'(D)\cdot\alpha^{-1}$. Similarly, for $s=1$, we end up with the relation $f(D)\dot p_1(\alpha)=1+\ln_{\mathfrak{T}f}\alpha$, which is indeed true, since we have an identity $f(D)\alpha=f'(D)+\alpha f(D)$ for operators and the definition of $\ln_f \alpha$.\\
There is one nice property, that the continuations of trivial polynomials for $f(x)=x \rightsquigarrow$ $\rightsquigarrow p_s(\alpha)=\alpha^s$ have:
$$\frac{\partial}{\partial\alpha} \frac{\alpha^s \ln\alpha}{\alpha^s}=\alpha^{-1}$$

Such an indentity is trivial by itself, but it should be noticed that the numerator is equal to the series $\dot p_s(\alpha)$, the denominator is equal to $p_s(\alpha)$, and the RHS of this identity is simply $p_{-1}(\alpha)$. There is no such an apparent structure in general, but in case of deformation $x \to \Delta_p=(e^{px}-1)/p$ such a property is preserved with minor changes, which follows from the properties of generalized harmonic numbers. But for the sake of correctness we need to avoid the use of identity
$$\Delta_p \rightsquigarrow p_s(\alpha)=p^s \frac{\Gamma(\frac{\alpha}{p}+1)}{\Gamma(\frac{\alpha}{p}+1-s)}$$
since we deal with the series, reflecting the asymptotic properties of $\Gamma$-function. In comparison, the series
$$p_s^*(\alpha)=(-p)^s \frac{\Gamma(-\frac{\alpha}{p}+s)}{\Gamma(-\frac{\alpha}{p})}$$
from our formal point of view has the same properties, but represents completely different class of functions, coinciding with the previous one only for integer values of $s$. Such a confusion appears when we deal with functions, not the series, because the kernel of operator $\Delta_p(D)$ is nontrivial and consists of all $p$-periodic functions. Although this distinction does not affect our study, the caution is still required. That is how the inaccurate use of the explicit formulae for general falling factorials $e^x-1 \rightsquigarrow (\alpha)_s$ for different \textit{functional} continuations leads to incoherent results with generalized harmonic numbers
$$
\frac{\partial}{\partial s}(\alpha)_s = (\alpha)_s (-\gamma+H_{\alpha-s}) ~~~~ \frac{\partial}{\partial s}(\alpha)_s = (\alpha)_s (\ln(-1)-\gamma+H_{s-\alpha-1}) 
$$
to make which agreeable one needs more than just fixing the meaningful value for logarithm of negative unit. It is still necessary to notice the difference between bracketed expressions, which in functional terms can be written as $\pi \cot (\pi (s-\alpha))$, the function that lies in the kernel of $e^D-1$ because of its 1-periodicity. Consider the following statement, in which we do not use the conventional functional notation on purpose.

We note only that the logarithm of any series of the form $\alpha^s+\alpha^{s-1}\mathbb{C}[[\alpha^{-1}]]$ is well-defined and is of the form $s\ln\alpha+A_1 \alpha^{-1}+...$ for the formal element $\ln\alpha$, as in general raising to any power $h$ is well-defined. Also it follows from the standard convolution identities that $e^{AD} g^h=(e^{AD}g)^h$ and, in particular, $e^{AD}\ln g=\ln(e^{AD}g)$ for any $g \in \alpha^s+\alpha^{s-1}\mathbb{C}[[\alpha^{-1}]]$.\\
~\\
\textbf{Proposition 2.2}. Suppose $f(x)=\Delta_p e^{-Ax} \rightsquigarrow p_s(\alpha)$. Then the following holds true:
\begin{align*}
\tag{2.5}f\left(\frac{\partial}{\partial\alpha}\right)\ln\frac{p_s(\alpha)}{\alpha}&=-f\left((1-s)\frac{\partial}{\partial\alpha}\right)\ln\alpha\\
\tag{2.6}f\left(\frac{\partial}{\partial\alpha}\right)\frac{\dot p_s(\alpha)}{p_s(\alpha)}&=\frac{1}{1-s}p_{-1}\left(\frac{\alpha}{1-s}\right)
\end{align*}
\emph{Proof}: Consider the case $A=0$, i.e. $\Delta_p(x) \rightsquigarrow p_s(\alpha)$. Since for $f(x) \rightsquigarrow p_s(\alpha)$, we have $\alpha (f'(D))^{-1}p_s(\alpha)=p_{s+1}(\alpha)$, in this case we have $\alpha e^{-pD} p_s(\alpha)=p_{s+1}(\alpha)$. Or, equivalently $e^{pD}p_{s+1}(\alpha)=(\alpha+p)p_{s}(\alpha)$. Then it follows that $(e^{pD}-1)p_{s+1}(\alpha)=(\alpha+p)p_{s}(\alpha)-p_{s+1}(\alpha)$. On the other hand, by definition, $(e^{pD}-1)p_{s+1}(\alpha)=p(s+1)p_s(\alpha)$, hence $p_{s+1}(\alpha)=$ $=(\alpha-ps)p_s(\alpha)$, as expected. Then
\begin{align*}
&\frac{e^{pD}-1}{p}\ln p_{s}(\alpha)=\frac{1}{p}\ln\left(\frac{e^{pD}p_s(\alpha)}{p_s(\alpha)}\right)=\frac{1}{p}\ln\left(\frac{(\alpha+p)p_{s-1}(\alpha)}{p_s(\alpha)}\right)=\\
&=\frac{1}{p}\ln\left(\frac{\alpha+p}{\alpha+p(1-s)}\right)=\frac{1}{p}\ln\left(\frac{e^{pD}\alpha}{e^{p(1-s)D}\alpha}\right)=\frac{e^{pD}-e^{p(1-s)D}}{p}\ln\alpha
\end{align*}
And so the following holds
$$
\Delta_p(D)\ln\frac{p_s(\alpha)}{\alpha}=\frac{1-e^{p(1-s)D}}{p}\ln\alpha=-\Delta_p((1-s)D)\ln\alpha
$$
It is left to say that the transformation $f(x) \to f(x)e^{-Ax}$, changes $p_s$ via the rule $p_s(\alpha) \to$ $\to \alpha e^{AsD}(p_s(\alpha)/\alpha)$, and thus the series
$$
\frac{f\left(\frac{\partial}{\partial\alpha}\right)}{f\left((1-s)\frac{\partial}{\partial\alpha}\right)}\ln\frac{p_{\vphantom{1} s}^f(\alpha)}{\alpha}
$$
remains unchanged under the transformation $f(x) \to f(x)e^{-Ax}$ for any $A$. So (2.5) indeed holds true. Now differentiating it by $s$ and taking into account the identity $f'(D)\cdot\alpha^{-1}=$ $=p_{-1}(\alpha)$, one can obtain
\begin{align*}
\tag*{\qed}&f\left(\frac{\partial}{\partial\alpha}\right)\frac{\dot p_s(\alpha)}{p_s(\alpha)}=f'\left((1-s)\frac{\partial}{\partial\alpha}\right)\frac{1}{\alpha}=\frac{1}{1-s}p_{-1}\left(\frac{\alpha}{1-s}\right)
\end{align*}
\textbf{Remark 2.4}. In case $s=1$ expression (2.6) still makes sense, because
$$
\frac{1}{1-s}p_{-1}\left(\frac{\alpha}{1-s}\right)=\alpha^{-1}+A_1(1-s)\alpha^{-2}+A_2(1-s)^2\alpha^{-3}+...
$$
for some $A_i$, hence it is equal to $\alpha^{-1}$ for $s=1$, and thus (2.6) follows from (2.4) in this case:
$$
f\left(\frac{\partial}{\partial\alpha}\right)\frac{\dot p_1(\alpha)}{p_1(\alpha)}=f\left(\frac{\partial}{\partial\alpha}\right)\ln_f \alpha=\frac{1}{\alpha}
$$
It should also be mentioned that (2.5), (2.6) for $s=0,1$ actually hold for any $f \in x+$ $+x^2\mathbb{C}[[x]]$. Although for other values of $s$ the general expansion $\ln(p_s(\alpha)/\alpha)=(s-1)\ln\alpha+$ $+A_1(s)\alpha^{-1}+...$ is difficult to study. Nevertheless, there exists at least one more regular case. Notice that $\mathfrak{T}^2 y_p=\mathfrak{T}^2 (\Delta_pe^{-x})=\Delta_p(1-\Delta_{-p})$ and consider the \mbox{following statement.}\\
~\\
\underline{\textbf{Theorem 2.1}}. Consider $f(x)=\mathfrak{T}^2(\Delta_p e^{-Ax})=\Delta_p(1-A\Delta_{-p}) \rightsquigarrow p_s(\alpha)$. Then the following holds:
\begin{align*}
\tag{2.7}f\left(\frac{\partial}{\partial\alpha}\right)\ln\frac{p_s(\alpha)}{p_{1-s}(\alpha)}=\left(f\left(s\frac{\partial}{\partial\alpha}\right)-f\left((1-s)\frac{\partial}{\partial\alpha}\right)\right)&\ln\alpha
\end{align*}
\begin{align*}
\tag{2.8}f\left(\frac{\partial}{\partial\alpha}\right)\left(\frac{\dot p_s(\alpha)}{p_s(\alpha)}+\frac{\dot p_{1-s}(\alpha)}{p_{1-s}(\alpha)}\right)=\frac{1}{s}p_{-1}\left(\frac{\alpha}{s\vphantom{1-s}}\right)+\frac{1}{1-s}p_{-1}&\left(\frac{\alpha}{\vphantom{s}1-s}\right)
\end{align*}
The proof is divided in a few parts.\\
~\\
\textbf{Proposition 2.3}. Consider $f(x)=\mathfrak{T}^2(\Delta_p e^{-Ax})=\Delta_p(1-A\Delta_{-p}) \rightsquigarrow p_s(\alpha)$, $\Delta_p \rightsquigarrow p_{\vphantom{1} s}^{\Delta_p}(\alpha)$. Then the following reflection formula holds true:
$$
p_{\vphantom{1} s}^{\Delta_p}(\alpha)p_{1-s}^{\vphantom{\Delta_p}}(\alpha)=p_{1-s}^{\Delta_p}(\alpha)p_{\vphantom{1} s}^{\vphantom{\Delta_p}}(\alpha)\eqno(2.9)
$$
\emph{Proof}: It is known that
$$e^{-pD}p_{\vphantom{1} s}^{\Delta_p}(\alpha)=\alpha^{-1}p_{s+1}^{\Delta_p}(\alpha) ~\Rightarrow~ e^{-npD}p_{\vphantom{1} s}^{\Delta_p}(\alpha)=\alpha^{-1}...(\alpha-(n-1)p)^{-1}p_{\vphantom{1} s+n}^{\Delta_p}(\alpha)$$
Thus
\begin{align*}
p_{\vphantom{1} s}(\alpha)&=\alpha \left(\frac{D}{f(D)}\right)^s \alpha^{s-1}=\alpha \left(\frac{\Delta_p(D)}{f(D)}\right)^s \left(\frac{D}{\Delta_p(D)}\right)^s\alpha^{s-1}=\\
&=\alpha (1-A\Delta_{-p})^{-s} \alpha^{-1}p_{\vphantom{1} s}^{\Delta_p}(\alpha)=\alpha\sum_{n=0}^\infty\binom{-s}{n}(-A)^n e^{-npD}\Delta_p^n(D)\alpha^{-1}p_{\vphantom{1} s}^{\Delta_p}(\alpha)
\end{align*}
Hence
\begin{align*}
p_s(\alpha)&=\alpha\sum_{n=0}^\infty\binom{-s}{n}(-A)^n e^{-npD}(s-1)_n\alpha^{-1}p_{\vphantom{1} s-n}^{\Delta_p}(\alpha)=\\
&=\alpha\sum_{n=0}^\infty\binom{s+n-1}{2n}\frac{(2n)!}{n!}A^n(\alpha-np)^{-1}\alpha^{-1}...(\alpha-(n-1)p)^{-1}p_{\vphantom{1} s}^{\Delta_p}(\alpha)=\\
&=p_{\vphantom{1} s}^{\Delta_p}(\alpha)\sum_{n=0}^\infty\binom{s+n-1}{2n}\frac{(2n)!}{n!}A^n(\alpha-p)^{-1}...(\alpha-np)^{-1}
\end{align*}
Then it follows that $p_s(\alpha)/p_s^{\Delta_p}(\alpha)$ is invariant under the transformation $s\to1-s$, since $((1-s)+n-1)...(1-s)((1-s)-1)...((1-s)-n)=(s+n-1)...s(s-1)...(s-n)$.\qed\\
~\\
\textbf{Remark 2.5}. In case $p=0$ this statement degenerates to the known property of $p_s$, associated to $x-x^2 \rightsquigarrow p_s(\alpha)=\frac{1}{\sqrt{\pi}}\alpha^{s+\frac{1}{2}}e^{\frac{\alpha}{2}}\mathrm{K}_{s-\frac{1}{2}}(\frac{\alpha}{2})$, where $\mathrm{K}_s(\alpha)$ denotes the modified Bessel function of the second kind (see \cite{Bes}):
$$
p_{1-s}(\alpha)=\alpha^{1-2s}p_{\vphantom{1} s}(\alpha)
$$
\textbf{Proposition 2.4}. $f(x)=\Delta_p(1-A\Delta_{-p}) ~\Rightarrow~ f(sx)-f((1-s)x)=(1-A\Delta_{-p})\cdot\vphantom{.}$ $\vphantom{.}\cdot(s\Delta_{ps}-(1-s)\Delta_{p(1-s)})$\\
\emph{Proof}: This property can be easily verified by expanding the brackets.\qed\\
~\\
Now use the latter propositions for $f(x)=\Delta_p(1-A\Delta_{-p})$ to obtain the following:
\begin{align*}
&f\left(\frac{\partial}{\partial\alpha}\right)\ln\frac{p_{\vphantom{1} s}(\alpha)}{p_{1-s}(\alpha)}=f\left(\frac{\partial}{\partial\alpha}\right)\ln\frac{p_{\vphantom{1} s}^{\Delta_p}(\alpha)}{p_{1-s}^{\Delta_p}(\alpha)}=\Delta_p(D)(1-A\Delta_{-p}(D))\ln\frac{p_{\vphantom{1} s}^{\Delta_p}(\alpha)}{p_{1-s}^{\Delta_p}(\alpha)}=\\
&=\frac{1-A\Delta_{-p}(D)}{p}\ln\frac{e^{pD}p_{\vphantom{1} s}^{\Delta_p}(\alpha)}{p_{\vphantom{1} s}^{\Delta_p}(\alpha)}\frac{p_{1-s}^{\Delta_p}(\alpha)}{e^{pD}p_{1-s}^{\Delta_p}(\alpha)}=\frac{1-A\Delta_{-p}(D)}{p}\ln\frac{\alpha+p}{\alpha+p(1-s)}\frac{\alpha+ps}{\alpha+p}=\\
&=\frac{e^{psD}-e^{p(1-s)D}}{p}(1-A\Delta_{-p}(D))\ln\alpha=\left(f\left(s\frac{\partial}{\partial\alpha}\right)-f\left((1-s)\frac{\partial}{\partial\alpha}\right)\right)\ln\alpha
\end{align*}
Take the derivative with respect to $s$ to obtain the remaining statement. \hfill\ensuremath{\blacksquare}\\
~\\
In this proof, the reflection formula $p_{\vphantom{1} s}^f(\alpha)p_{1-s}^g(\alpha)=p_{\vphantom{1} s}^g(\alpha) p_{1-s}^f(\alpha)$ for $g=\Delta_p(1-A\Delta_{-p})$, $f=\Delta_p$ plays a significant role. We do not consider the general equation $p_{\vphantom{1} s}^{A_1}(\alpha)p_{1-s}^{A_2}(\alpha)=$ $=p_{\vphantom{1} s}^{A_3}(\alpha)p_{1-s}^{A_4}(\alpha)$, one set of solutions of which includes $A_1=\Delta_p, A_2=\Delta_{-p}, A_3=\Delta_q$, $A_4=\Delta_{-q}$ (hence $p_{\vphantom{1} s}^{\Delta_p}(\alpha)p_{1-s}^{\Delta_{-p}}(\alpha)=\alpha$), in this paper, but the solution of special case $p_{\vphantom{1} s}^f(\alpha)p_{1-s}^g(\alpha)=p_{\vphantom{1} s}^g(\alpha) p_{1-s}^f(\alpha)$ is provided by the following theorem.\\
~\\
\underline{\textbf{Theorem 2.2}}. Suppose $f\neq g$ and $p_{\vphantom{1} s}^f(\alpha)p_{1-s}^g(\alpha)=p_{\vphantom{1} s}^g(\alpha) p_{1-s}^f(\alpha)$. Then there exist $p, A, B$, such that $f(x)=\Delta_p(1+A\Delta_{-p})$, $g(x)=\Delta_p(1+B\Delta_{-p})$, or, in other words, this equation does not have any solutions other than those considered above.\\
The proof is divided in a few parts.\\
~\\
\textbf{Proposition 2.5}.
$$
p_{\vphantom{1} s}^f(\alpha)p_{1-s}^g(\alpha)=p_{\vphantom{1} s}^g(\alpha) p_{1-s}^f(\alpha) ~\Rightarrow~ g(x)=f(x)(1+A\mathfrak{T}^{-1}f(x))
$$
for some constant $A$, or, equivalently $g(x)=\mathfrak{T}^2 e^{Ax} \mathfrak{T}^{-2}f(x)$.
\newpage
\noindent\emph{Proof}: Take the logarithmic derivative of LHS and RHS of this equation and set $s=0$ (or, equivalently $s=1$). According to Proposition 2.1, we have
$$
\dot p_0^f-\dot p_1^g\alpha^{-1}=\dot p_0^g-\dot p_1^f\alpha^{-1}  \iff \ln_{\mathfrak{T}f}\alpha - \ln_{g}\alpha = \ln_{\mathfrak{T}g}\alpha -\ln_{f}\alpha
$$
Hence the following holds
$$
\frac{xf'(x)}{f(x)}-\frac{x}{g(x)}=\frac{xg'(x)}{g(x)}-\frac{x}{f(x)}
$$
by definition of $\ln_M \alpha$ for any $M(x)$. Now integrate this expression, after dividing by $x$:
$$
C+\ln(f(x)\cdot \mathfrak{T}^{-1} f(x))=\ln(g(x)\cdot \mathfrak{T}^{-1}g(x))
$$
But $f,g,\mathfrak{T}^{-1} f, \mathfrak{T}^{-1} g \in x+x^2\mathbb{C}[[x]]$, and thus $C=0$. Hence
$$
f(x)\cdot \mathfrak{T}^{-1} f(x)=g(x)\cdot \mathfrak{T}^{-1}g(x)
$$
By definition, $f=\mathfrak{T}^{-1} f/(\mathfrak{T}^{-1} f)'$, so we can rewrite the expression above again
$$
\frac{(\mathfrak{T}^{-1} g(x))'}{(\mathfrak{T}^{-1} g(x))^2}=\frac{(\mathfrak{T}^{-1} f(x))'}{(\mathfrak{T}^{-1} f(x))^2} \iff \frac{1}{\mathfrak{T}^{-1} g(x)}= A + \frac{1}{\mathfrak{T}^{-1} f(x)}
$$
Hence $\mathfrak{T}^{-2}g(x)= e^{Ax} \mathfrak{T}^{-2}f(x)$ or, equivalently, $g(x)=f(x)(1+A\mathfrak{T}^{-1}f(x))$\qed\\
~\\
\textbf{Proposition 2.6}.
$$
p_{\vphantom{1} s}^f(\alpha)p_{1-s}^g(\alpha)=p_{\vphantom{1} s}^g(\alpha) p_{1-s}^f(\alpha) ~\Rightarrow~ f'(x)+g''(0)f(x)=g'(x)+f''(0)g(x)
$$
\emph{Proof}: Substituting the value $s=2$ (or, equivalently $s=-1$) into the equation and using the simple fact that $p_2^f(\alpha)=\alpha(\alpha-f''(0))$, we obtain
$$
\alpha(\alpha-f''(0))p_{-1}^g(\alpha)=\alpha(\alpha-g''(0))p_{-1}^f(\alpha)
$$
Taking into account the formulae $p_{-1}^f(\alpha)=-\alpha f(D) \cdot \alpha^{-1}$ and $\alpha g(D)=-g'(D)+g(D)\alpha$, we can rewrite the expression above in the following manner:
$$
(-g'(D)+g(D)(\alpha-f''(0)))\frac{1}{\alpha}=(-f'(D)+f(D)(\alpha-g''(0)))\frac{1}{\alpha}
$$
Also $f(D)1=g(D)1=0$, thus
$$
(g'(D)+f''(0)g(D))\frac{1}{\alpha}=(f'(D)+g''(0)f(D))\frac{1}{\alpha}
$$
and hence $f'(x)+g''(0)f(x)=g'(x)+f''(0)g(x)$\qed\\
~\\
Now we are ready to prove the theorem. Substituting $g=f(1+A\mathfrak{T}^{-1}f)$, since $(\mathfrak{T}^{-1}f)'=$ $=(\mathfrak{T}^{-1}f)/f$ and $g''(0)=f''(0)+2A$, one can easily obtain:
$$
f'(x)+f(x)(f''(0)+2A)=[f(x)(1+A\mathfrak{T}^{-1}f(x))]' + f''(0)f(x)(1+A\mathfrak{T}^{-1}f(x)) \Rightarrow
$$
$$
\Rightarrow~2Af(x)=Af'(x)\mathfrak{T}^{-1}f(x)+A\mathfrak{T}^{-1}f(x)+Af''(0)f(x)\mathfrak{T}^{-1}f(x)
$$
The latter can be rewritten in the following form:
$$
2A\frac{1}{\mathfrak{T}^{-1}f}(x)=A\frac{f'}{f}(x)+A\frac{1}{f}(x)+Af''(0)
$$
\newpage
\noindent Now, since we are interested in nontrivial solutions, we can divide both sides by $A\neq 0$ and integrate. Then it follows that
$$
\frac{f\cdot\mathfrak{T}^{-1}f}{(\mathfrak{T}^{-2}f)^2}=e^{-f''(0)x}
$$
But such an equation is equivalent to the following one:
$$
(q')^2-qq''=e^{f''(0)x} ~~\text{for}~~ q=\mathfrak{T}^{-2}f
$$
Differentiating and rearranging the summands, we obtain:
$$
q'(q''-f''(0)q')=q(q'''-f''(0)q'')
$$
Now if $q''-f''(0)q'=0$, then $q(x)=\Delta_{f''(0)}(x)$. Otherwise, in case $q''-f''(0)q'\neq0$, there exists some constant $C$, such that
$$
\frac{q'}{q}=\frac{q'''-f''(0)q''}{q''-f''(0)q'} ~\Rightarrow~ Cq(x)=q''(x)-f''(0)q'(x)
$$
Then it follows that $\exists c_1, c_2: q(x)=(e^{c_1 x}-1)e^{-c_2 x}/c_1$, and that includes the case $q''=$ $=f''(0)q'$. Hence $f$ can be expressed as
$$
f(x)=\frac{e^{c_1 x}-1}{c_1}\left(1-c_2\frac{1-e^{-c_1 x}}{c_1}\right)
$$
Thus $g$ can be expressed as
$$
g(x)=\frac{e^{c_1 x}-1}{c_1}\left(1+(A-c_2)\frac{1-e^{-c_1 x}}{c_1}\right)
$$
and that proves the theorem.\hfill\ensuremath{\blacksquare}


\begin{center}
\textbf{General comment on section 2}
\end{center}
The formula for $\dot p_s/p_s$, successfully found in case $y_p \rightsquigarrow p_s$ is not surprising, although it is not obvious, if it has any natural generalization. In case $f(x)=e^x-1 \rightsquigarrow p_n(\alpha)$ there exists an identity $(\alpha e^{-D})^n=p_n(\alpha) e^{-nD}$ for operators, and thus the operator \mbox{$(\alpha\frac{1}{f'}(D))^n=(\alpha e^{-D})^n$} is proportional to the operator of multiplication by $p_n$. Such a property is important, because in this case the operator of division by $p_n(\alpha)$ can be written down explicitly like $e^{-nD}(\alpha e^{-D})^{-n}$. Taking into account that for all $s$ we have \mbox{$(\alpha e^{-D})^n p_s(\alpha)=p_{s+n}(\alpha)$}, we obtain that also $(\alpha e^{-D})^{-n} \dot p_n(\alpha)=\dot p_0(\alpha) ~\Rightarrow~ \dot p_n(\alpha)/p_n(\alpha)=\dot p_0(\alpha-n)$. That is why we can write down the expression for $f(D)\dot p_n(\alpha)/p_n(\alpha)$ in this case. On the other hand, in case $f=\mathfrak{T}^2 y_p$ there is no such apparent relation between the shift operator $A_f^n\coloneqq(\alpha\frac{1}{f'}(D))^n$ and corresponding polynomials, but there is another equality of operators
$$
A_f^{n+1} f'(D)^{2n+1} A_f^{n\vphantom{+1}}=p^{2n+1}\left(\frac{\alpha}{p}+n\right)_{2n+1}\eqno{(2.10)}
$$
for falling factorials. We can rewrite the latter as
$$
p^{2n+1}\left(\frac{\alpha}{p}+n\right)_{2n+1} A_f^{-n}=A_f^{n+1} f'(D)^{2n+1}
$$
and act with both sides on $1$ to obtain the familiar formula
$$
p_{\vphantom{1} n}^{\Delta_p}(\alpha)p_{1-n}^{\vphantom{\Delta_p}}(\alpha)=p_{1-n}^{\Delta_p}(\alpha)p_{\vphantom{1} n}^{\vphantom{\Delta_p}}(\alpha)
$$
An interested reader might find the proof of (2.10) in Appendix B.

\newpage
\hrule

\section{On the basic properties of derivation and shift operators}

\hrule

~\\
~\\
In this section we study the main properties of operators $\alpha f'(D)^{-1}$ and $f(D)$, which we call \emph{shift} and \emph{derivation} respectively. The properties provided allow us to understand the behaviour of polynomials $p_n^g$ under the transformations $\mathfrak{T}^{-1}f \to f \to \mathfrak{T}f$ of associated series more fully. In conclusion, we introduce two additional formulae for the sequence $a_n$, considered in the previous paper.\\
As usual we denote the falling factorials $x(x-1)...(x-n+1)$ by $(x)_n$. For conciseness, the operators $\alpha f'(D)^{-1}$ and $f(D)$ are sometimes denoted by $A_f$ and $D_f$ respectively.\\
~\\
\textbf{Proposition 3.1}. Consider some polynomial $q_n(\alpha)$ of degree $n$, which has the coefficient of leading term equal to $1$. Suppose $\mathfrak{T}^{-1}f(x) \rightsquigarrow p_n^{\mathfrak{T}^{-1}f}(\alpha)$. Then the following holds:
$$
p_n^{\mathfrak{T}^{-1}f}(\alpha)=\frac{1}{n!}\left(\alpha f(D)\right)_n \cdot q_n(\alpha)\eqno(3.1)
$$
In particular,
$$
p_n^{\mathfrak{T}^{-1}f}(\alpha)=\frac{1}{n!}\left(\alpha f(D)\right)_n \cdot p_n^f(\alpha)
$$
\emph{Proof}: The operators $A_f$, $D_f$ clearly satisfy $D_f A_f =1+A_f D_f$. Then it folows from this relation that $A_f^{n\vphantom{-1}} D_f^{n\vphantom{-1}}=\left(A_f D_f\right)_n$. To be precise: the latter obviously holds in case $n=1$, and if it holds for some $k$, from the simple corollary
$$D_f A_f =1+A_f D_f \Rightarrow D_f^{n\vphantom{-1}} A_f^{\vphantom{n}}=nD_f^{n-1}+A_f^{\vphantom{n}}D_f^{n\vphantom{-1}}$$
it follows that
\begin{align*}
\left(A_f D_f\right)_{k+1}&=\left(A_f D_f\right)_k (A_f D_f -k)=A_f^{k\vphantom{-1}} D_f^{k\vphantom{-1}}A_f D_f-kA_f^{k\vphantom{-1}} D_f^{k\vphantom{-1}}=\\
&=A_f^{k\vphantom{-1}}(kD_f^{k-1}+A_f^{\vphantom{-1}} D_f^{k\vphantom{-1}})-kA_f^{k\vphantom{-1}}D_f^{k\vphantom{-1}}=A_f^{k+1}D_f^{k+1}
\end{align*}
Hence, indeed $A_{\mathfrak{T}^{-1}f}^{n\vphantom{-1}} D_{\mathfrak{T}^{-1}f}^{n\vphantom{-1}} = \left(\alpha f(D)\right)_n$. Acting with the latter operator on the polynomial $q_n(\alpha)$, we obtain: $\left(\alpha f(D)\right)_n \cdot q_n(\alpha)=A_{\mathfrak{T}^{-1}f}^{n\vphantom{-1}} D_{\mathfrak{T}^{-1}f}^{n\vphantom{-1}}\cdot q_n(\alpha)=n!A_{\mathfrak{T}^{-1}f}^{n\vphantom{-1}}\cdot 1=n!p_n^{\mathfrak{T}^{-1}f}(\alpha)$. In particular, setting $q_n(\alpha)=p_n^f(\alpha)$, we obtain the remaining statement.\qed\\
~\\
The latter proposition shows the way to express polynomials $p_n^{\mathfrak{T}^{-1}f}(\alpha)$, which does not require the Lagrange inversion theorem and computation of the series of powers $(\mathfrak{T}^{-1}f)^{-k}$. We use this method below to provide some new formulae for the sequence $a_n$:
$$
\sum_{n=1}^\infty a_n x^n=\left(\mathfrak{T}^{-1} xe^{-x}\right)^{inv}
$$
But for now consider a few more statements.\\
~\\
\textbf{Proposition 3.2}. (\emph{formal}). Suppose that the action of operator $f'(D)^{-1}$ is well-defined not only in case of the series $\in \alpha^s \mathbb{C}[[\alpha^{-1}]]$, but in more general case, as for polynomial $f'(x)^{-1}=P(x)$ or exponent $f'(x)^{-1}=e^{px}$. Suppose $f(x)$ has a simple zero $\rho$: $f(\rho)=0$, $f'(\rho)\neq 0$ and has an expansion $f(x+\rho)\approx f'(\rho)x+A_1 x^2+...$ at that point. Consider the polynomials $f(x+\rho)(f'(\rho))^{-1}\rightsquigarrow p_{s\vphantom{-1}}^{\rho}(\alpha)$. Then if $A_f=\alpha f'(D)^{-1}$, the following formal identity holds true:
$$
p_{n\vphantom{-1}}^{\rho}(\alpha)=f'(\rho)^n e^{-\alpha\rho}A_f^n \cdot e^{\alpha\rho}\eqno(3.2)
$$
\newpage
\noindent\emph{Proof}: If the action of $f'(D)^{-1}$ on exponent is well-defined, then the result must be of the following canonical form: $f'(D)^{-1}\cdot e^{\alpha t}=f'(t)^{-1}e^{\alpha t}$. Expanding the operator $A_f^n$ in sum $\sum_{k=0}^n \alpha^k g_k^n(D)$ where $g$ are some functions, depending on $f$, one can obtain:
$$
f'(\rho)^n e^{-\alpha\rho}A_f^n \cdot e^{\alpha\rho}=f'(\rho)^n\sum_{k=0}^n \alpha^k g_k^n(\rho)
$$
and, moreover, $g_k^n(\rho)$ depend only on coefficients of the expansion of $f$ at point $\rho$. On the other hand such an expression is nothing more than the expansion of
$$
f'(\rho)^n\sum_{k=0}^n \alpha^k g_k^n(\rho)=f'(\rho)^n\left(\alpha\frac{1}{f'(\rho+D)}\right)^n \cdot 1=\left(\alpha\frac{f'(\rho)}{f'(\rho+D)}\right)^n \cdot 1
$$
and that is indeed $p_{n\vphantom{-1}}^{\rho}(\alpha)$.\qed\\
~\\
\textbf{Remark 3.1}. Of course, the latter proposition makes sense only in special cases, for example, if $f(x)=e^x-1$. In such ``nice'' cases the following equality of operators holds
$$
f\left(\rho+\frac{\partial}{\partial\alpha}\right)e^{-\alpha\rho}=e^{-\alpha\rho}f\left(\frac{\partial}{\partial\alpha}\right)
$$
which of course does not hold, when the action on exponent is not well-defined.\\
~\\
\textbf{Proposition 3.3}. Suppose $f\rightsquigarrow p_n^f$, $\mathfrak{T}f\rightsquigarrow p_n^{\mathfrak{T}f}$, $\mathfrak{T}^{-1}f\rightsquigarrow p_n^{\mathfrak{T}^{-1}f}$.
Then the following equalities of operators hold:
\begin{align*}
A_f^n f'(D)^n=\alpha\left(\alpha-\frac{f''(D)}{f'(D)}\right)\left(\alpha-2\frac{f''(D)}{f'(D)}\right)...\left(\alpha-(n-1)\frac{f''(D)}{f'(D)}\right)
\end{align*}
\begin{align*}
f'(D)^n A_f^n=\left(\alpha+n\frac{f''(D)}{f'(D)}\right)\left(\alpha+(n-1)\frac{f''(D)}{f'(D)}\right)...\left(\alpha+\frac{f''(D)}{f'(D)}\right)
\end{align*}
and the following expressions are valid:
\begin{align*}
\tag{3.3}\frac{p_{n+1}^f(\alpha)}{\alpha}&=\left(\alpha-\frac{f''(D)}{f'(D)}\right)...\left(\alpha-n\frac{f''(D)}{f'(D)}\right)\cdot 1\\
\tag{3.4}\frac{p_{n+1}^{\mathfrak{T}f}(\alpha)}{\alpha}&=\left(\alpha+n\frac{f''(D)}{f'(D)}\right)...\left(\alpha+\frac{f''(D)}{f'(D)}\right)\cdot 1\\
\tag{3.5}\frac{p_{n+1}^{\mathfrak{T}^{-1}f}(\alpha)}{\alpha}&=\left(\alpha+\frac{f'(D)-1}{f(D)}\right)...\left(\alpha+n\frac{f'(D)-1}{f(D)}\right)\cdot 1
\end{align*}
(as one may notice, in this case the operators do not commute in general)\\
\emph{Proof}: It should be mentioned that $f'(D)^np_{n\vphantom{+1}}^f(\alpha)=\alpha^{-1}p_{n+1}^{\mathfrak{T}f}(\alpha)$. Now consider the following equality of operators:
$$
f'(D)^\lambda\alpha f'(D)^{-\lambda}=\alpha+\lambda\frac{f''(D)}{f'(D)}
$$
Then
\begin{align*}
\alpha\left(\alpha-\frac{f''(D)}{f'(D)}\right)&\left(\alpha-2\frac{f''(D)}{f'(D)}\right)...\left(\alpha-(n-1)\frac{f''(D)}{f'(D)}\right)=\\
&=\alpha f'(D)^{-1}\alpha f'(D)\cdot f'(D)^{-2}\alpha f'(D)^2\cdot...\cdot f'(D)^{-n+1}\alpha f'(D)^{n-1}=\\
&=\left(\alpha\frac{1}{f'(D)}\right)^n f'(D)^n=A_f^n f'(D)^n
\end{align*}
\begin{align*}
\left(\alpha+n\frac{f''(D)}{f'(D)}\right)&\left(\alpha+(n-1)\frac{f''(D)}{f'(D)}\right)...\left(\alpha+\frac{f''(D)}{f'(D)}\right)=\\
&=f'(D)^{n}\alpha f'(D)^{-n}\cdot f'(D)^{n-1}\alpha f'(D)^{1-n}\cdot...\cdot f'(D)\alpha f'(D)^{-1} =\\
&=f'(D)^n\left(\alpha\frac{1}{f'(D)}\right)^n=f'(D)^nA_f^n
\end{align*}
Now act with the first operator on $1$. Since $f'(D)^n\cdot 1=1$, we obtain
$$
\alpha\left(\alpha-\frac{f''(D)}{f'(D)}\right)\left(\alpha-2\frac{f''(D)}{f'(D)}\right)...\left(\alpha-(n-1)\frac{f''(D)}{f'(D)}\right)\cdot 1=p_{n\vphantom{+1}}^f(\alpha)
$$
Then after change of index $n \to n+1$ we end up with (3.3). Similarly,
\begin{align*}
\left(\alpha+n\frac{f''(D)}{f'(D)}\right)\left(\alpha+(n-1)\frac{f''(D)}{f'(D)}\right)...\left(\alpha+\frac{f''(D)}{f'(D)}\right)\cdot 1=f'(D)^np_{n\vphantom{+1}}^f(\alpha)=\alpha^{-1}p_{n+1}^{\mathfrak{T}f}(\alpha)
\end{align*}
Now (3.5) follows from (3.3), since we have, by definition,
\begin{align*}
\tag*{\qed}\frac{(\mathfrak{T}^{-1}f)''}{(\mathfrak{T}^{-1}f)'}=\frac{1-f'}{f}
\end{align*}
On this note we finish the overview of basic identities and return to Proposition 3.1. Consider the function $x e^{\mathrm{Ei}(x)-\ln|x|-\gamma}\rightsquigarrow \nu_n(\alpha)$, which was studied in the previous paper, and the binomial polynomials, associated to this function. We remind the reader, that these polynomials have the property $n\nu_n(\alpha)=\alpha \nu_n'(\alpha-1)$ and have the following explicit representation:
$$
\nu_n(\alpha)=\sum_{k=0}^{n-1}\binom{n-1}{k}\alpha^{n-k}\sum_{m=0}^{k}\frac{(-n)^m}{m!}\sum_{\substack{\sum_{1}^{m}q_i=k\\ q_i>0}}\binom{k}{q_1,...,q_k}\frac{1}{q_1...q_k}
$$
In particular, we are interested in the numbers $a_n=\nu_n'(0)/n!$. As one may notice, their generating function is the series $\psi(x)$, considered in Proposition 1.9.\\
~\\
\underline{\textbf{Theorem 3.1}}. Suppose $xe^{\mathrm{Ei}(x)-\ln|x|-\gamma}\rightsquigarrow \nu_n(\alpha)$, $a_n=\nu_n'(0)/n!$, and $\begin{bsmallmatrix} n\\m \end{bsmallmatrix}$ are the signed Stirling numbers of the first kind (i.e. the coefficients $[x^m](x)_n$). Then for any $n\geqslant 1$ the following holds
\begin{align*}
&(-1)^{n-1}(n-1)!\nu_n(\alpha)\alpha^{-1}=\\
\tag{3.6}&=\sum_{k=1}^n \begin{bmatrix} \phantom{i}n\phantom{i}\\ \phantom{i}k\phantom{i} \end{bmatrix}\smashoperator[r]{\sum_{\sum_1^k m_i=n-1}} \binom{n-1}{m_1,...,m_k}(n-m_1)...(n-m_1-...-m_{k-1})(1-\alpha)^{m_k}\\
\tag{3.7}&(-1)^{n-1}n\nu_n(\alpha)\alpha^{-1}=\sum_{k=1}^n \begin{bmatrix} \phantom{i}n\phantom{i}\\ \phantom{i}k\phantom{i} \end{bmatrix}\smashoperator[r]{\sum_{\sum_1^{k+1} m_i=n-1}} \frac{n}{m_1+1}\frac{n-m_1}{m_2+1}...\frac{n-m_1-...-m_{k-1}}{m_k+1}\binom{m_{k+1}-\alpha}{m_{k+1}}
\end{align*}
In particular,
\begin{align*}
\tag{3.8}&(-1)^{n-1}n!(n-1)!a_n=\sum_{k=1}^n \begin{bmatrix} \phantom{i}n\phantom{i}\\ \phantom{i}k\phantom{i} \end{bmatrix}\smashoperator[r]{\sum_{\sum_1^k m_i=n-1}} \binom{n-1}{m_1,...,m_k}(n-m_1)...(n-m_1-...-m_{k-1})\\
\tag{3.9}&(-1)^{n-1}nn!a_n=\sum_{k=1}^n \begin{bmatrix} \phantom{i}n\phantom{i}\\ \phantom{i}k\phantom{i} \end{bmatrix}\smashoperator[r]{\sum_{\sum_1^{k+1} m_i=n-1}} \frac{n}{m_1+1}\frac{n-m_1}{m_2+1}...\frac{n-m_1-...-m_{k-1}}{m_k+1}
\end{align*}
\emph{Proof}: Since $\mathfrak{T}~x e^{\mathrm{Ei}(x)-\ln|x|-\gamma}=xe^{-x}$, according to (3.1), for any polynomial $q_n(\alpha)$ of degree $n$ with the coefficient of leading term equal to $1$, we have $n!\nu_n(\alpha)=(\alpha De^{-D})_n q_n(\alpha)$. Consider two cases: $q_n(\alpha)=\alpha^n$ and $q_n(\alpha)=(\alpha)_n$. The main idea of this proof is to rewrite the action of $\left(\alpha De^{-D}\right)^k$ in the following manner:
$$
\left.\frac{\partial}{\partial A_1}...\frac{\partial}{\partial A_k}\alpha e^{A_1 D} ... \alpha e^{A_k D}\right|_{A_i=-1}
$$
And this expression can be rewritten again as
$$
\left.\frac{\partial}{\partial A_1}...\frac{\partial}{\partial A_k}\alpha(\alpha+A_1)...(\alpha+A_1+...+A_{k-1}) e^{(A_1+...+A_k) D}\right|_{A_i=-1}
$$
and after action with this operator on $q_n(\alpha)$ we obtain:
\begin{align*}
&n!\nu_n(\alpha)=\\
&=\sum_{k=1}^n \begin{bmatrix} \phantom{i}n\phantom{i}\\ \phantom{i}k\phantom{i} \end{bmatrix}\left.\frac{\partial^k}{\partial A_1...\partial A_k}\prod_{j=0}^{k-1}\left(\alpha+\sum_{i=1}^{j}A_i\right)q_n(\alpha+A_1+...+A_k)\right|_{A_i=-1}
\end{align*}
Now use the binomiality of $q_n(\alpha)$. Consider the case $q_n(\alpha)=\alpha^n$. Then
$$
q_n(\alpha+A_1+...+A_k)=\sum_{m_1=0}^{n}\binom{n}{m_1}A_k^{m_1}q_{n-m_1}(\alpha+A_1+...+A_{k-1})
$$
That means after multiplication by $(\alpha+A_1+...+A_{k-1})$ we end up with the expression:
$$
\sum_{m_1=0}^{n}\binom{n}{m_1}A_k^{m_1}q_{n-m_1+1}(\alpha+A_1+...+A_{k-1})
$$
in which we may now change index $m_1 \to m_1+1$. Expanding in the same way again for $A_{k-1}$, after multiplication by $(\alpha+A_1+...+A_{k-2})$, and using the same strategy for next $A_i$, we obtain the expression
\begin{align*}
&\prod_{j=0}^{k-1}\left(\alpha+\sum_{i=1}^{j}A_i\right)q_n(\alpha+A_1+...+A_k)=\smashoperator{\sum_{\substack{-1\leqslant m_1 \leqslant n-1 \\ -1\leqslant m_2 \leqslant n-1-m_1\\...\\-1 \leqslant m_{k-1}\leqslant n-1-m_1-...-m_{k-2}}}} \mathrm{G}_{\mathrm{m}_i}^n
\end{align*}
where $\mathrm{G}_{\mathrm{m}_i}^n$ are the products of coefficients
$$
\binom{n}{m_1+1}\binom{n-m_1}{m_2+1}...\binom{n-m_1-...-m_{k-2}}{m_{k-1}+1}
$$
and polynomials
$$
A_k^{m_1+1}...A_2^{m_{k-1}+1}\alpha(\alpha+A_1)^{n-m_1-...-m_{k-1}}
$$
Now after differentiating by $A_i$, the limits of summation change from $-1\leqslant m_j$ to $0 \leqslant m_j$. Also each summand multiplies by $(1+m_1)...(1+m_{k-1})(n-m_1-...-m_{k-1})$. Then after substitution $A_j=-1$ we obtain
\begin{align*}
\left.\frac{\partial^k}{\partial A_1...\partial A_k}\prod_{j=0}^{k-1}\left(\alpha+\sum_{i=1}^{j}A_i\right)q_n(\alpha+A_1+...+A_k)\right|_{A_i=-1}=\smashoperator[r]{\sum_{\substack{0\leqslant m_1 \leqslant n-1 \\ 0\leqslant m_2 \leqslant n-1-m_1\\...\\0 \leqslant m_{k-1}\leqslant n-1-m_1-...-m_{k-2}}}} \mathrm{R}_{\mathrm{m}_i}^n
\end{align*}
where now $\mathrm{R}_{\mathrm{m}_i}^n$ are the products of the same coefficients
$$
\binom{n}{m_1+1}\binom{n-m_1}{m_2+1}...\binom{n-m_1-...-m_{k-2}}{m_{k-1}+1}
$$
and polynomials
$$
(1+m_1)...(1+m_{k-1})(n-m_1-...-m_{k-1})(-1)^{n-1}\alpha(1-\alpha)^{n-1-m_1+...+m_{k-1}}
$$
Now put the new variable $m_k=n-1-m_1+...+m_{k-1}$ to obtain the expression
$$
(-1)^{n-1}\smashoperator{\sum_{\sum_1^k m_i=n-1}}\binom{n}{m_1+1}\binom{n-m_1}{m_2+1}...\binom{n-m_1-...-m_{k-2}}{m_{k-1}+1}\alpha(1-\alpha)^{m_k}\prod_{i=1}^{k}(1+m_i)
$$
Since $\binom{n}{k}=\frac{n}{k}\binom{n-1}{k-1}$, we can rewrite the summands again as following
$$
(-1)^{n-1}\smashoperator{\sum_{\sum_1^k m_i=n-1}}\binom{n-1}{m_1,...,m_k}n(n-m_1)...(n-m_1-...-m_{k-1})\alpha(1-\alpha)^{m_k}
$$
Multiply each of this expression by the corresponding Stirling number and take the sum to obtain (3.6).\\
Consider the second case $q_n(\alpha)=(\alpha)_n$. We now use the same strategy, but with minor changes. In this case we have $q_n(\alpha-1)\cdot\alpha=q_{n+1}(\alpha)$ and thus one should expand the summands for $A_k$ with shift:
\begin{align*}
&(\alpha+A_1+...+A_{k-1})(\alpha+A_1+...+A_k)_n=\\
&=\sum_{k=0}^{n}\binom{n}{m_1}(A_k+1)_{m_1}(\alpha+A_1+...+A_{k-1})_{n-m_1+1}
\end{align*}
Using this method for each new $A_i$, we end up with the expressions
\begin{align*}
\alpha(\alpha+A_1)...(\alpha+A_1+...+A_{k-1})(\alpha+A_1+...+A_k)_n~~=~~\smashoperator{\sum_{\substack{-1\leqslant m_1 \leqslant n-1 \\ -1\leqslant m_2 \leqslant n-1-m_1\\...\\-1 \leqslant m_{k}\leqslant n-1-m_1-...-m_{k-1}}}} \mathrm{G}_{\mathrm{m}_i}^n
\end{align*}
where now $\mathrm{G}_{\mathrm{m}_i}^n$ are the products of coefficients
$$
\binom{n}{m_1+1}\binom{n-m_1}{m_2+1}...\binom{n-m_1-...-m_{k-1}}{m_{k}+1}
$$
and polynomials
$$
(A_k+1)_{m_1+1}...(A_1+1)_{m_k+1}(\alpha)_{n-m_1-...-m_k}
$$
(note: the number of indices has changed, since we now expand the brackets for $A_1$, in constrast to the previous case)\\
Now differentiate by $A_i$, taking into account that $\frac{d}{dx}(x)_n|_{x=0}=(-1)^{n-1}(n-1)!$. As in previous case, the limits of summation change. That means after excluding the trivial summands with indices $m_i=-1$ and reducing the factorials, we obtain the following polynomials for new variable $m_{k+1}=n-1-m_1+...+m_{k}$:
$$
(-1)^{n-1}(n-1)!\alpha\smashoperator{\sum_{\sum_1^{k+1} m_i=n-1}}\frac{n(n-m_1)...(n-m_1-...-m_{k-1})}{(m_1+1)(m_2+1)...(m_k+1)}(-1)^{m_{k+1}}\binom{\alpha-1}{m_{k+1}}
$$
Now multiply each of these polynomial by the corresponding Stirling number and take the sum to obtain (3.7). The identities (3.8) and (3.9) then are the direct corollaries of (3.6) and (3.7) in special case $\alpha=0$.\hfill\ensuremath{\blacksquare}

\begin{center}
\textbf{General comment on section 3}
\end{center}
Going back to Proposition 3.2, in which we considered the formal identity for the polynomials, associated to the expansion of the same function at other points, it should be noticed that one can obtain more precise statement. For that purpose one can look at more general structure of formal series $\alpha^s e^{\alpha t} \mathbb{C}[[\alpha^{-1}]]$ for some formal element $t$. Now one can \emph{define} the action of $f(D)$ on exponent, considering $f(x)$ as an element of more general structure $\mathbb{C}[x,t]$. Then we may consider the following expansion by the variable $x$ with formal power series of $t$ as coefficients:
$$
\left(\frac{x}{f(x+t)-f(t)}\right)^s=\sum_{n=0}^\infty \frac{q_n^{t}(s)x^n}{n!}
$$
and then for any complex $s$ define
$$
A_f^s\coloneqq\sum_{n=0}^\infty \binom{s-1}{n}\alpha^{s-n}q_n^D(s)
$$
Now one can verify the following operator identities:
$$
D_f^{\vphantom{s-1}} A_f^{s\vphantom{-1}}=sA_f^{s-1}+A_f^{s\vphantom{-1}}D_f^{\vphantom{s-1}}; ~~A_f^{s\vphantom{+1}} A_f^{h\vphantom{s+1}} =  A_f^{h\vphantom{s+1}}A_f^{s\vphantom{+1}} = A_f^{s+h\vphantom{+1}}
$$
and the evaluation
$$
A_f^{s\vphantom{+1}}\cdot e^{\alpha t}=f'(t)^{-s}p_s^{t}(\alpha)e^{\alpha t}
$$
where $p_s^{t}(\alpha)$ for $t=0$ are equal to the ordinary canonical continuations of polynomials, and in general are naturally given by definition
$$
p_s^{t}(\alpha)\coloneqq\sum_{n=0}^\infty\binom{s-1}{n}\alpha^{s-n}q_n^t(s)f'(t)^s
$$
In particular, such an approach allows us to define the operator $\ln A_f$, for which the following holds:
$$
\ln A_f \cdot p_s(\alpha)=\dot p_s(\alpha)
$$
With regard to Theorem 3.1, one may derive, actually, more general statement from Proposition 3.1, which can be used as an alternative to the Lagrange inversion theorem, precisely:
$$
\mathrm{Q}\mathfrak{T}^{-1}\sum_{n=1}^\infty c_{n} x^{n}=\sum_{n=1}^\infty \frac{x^n}{n!} \sum_{k=1}^{n}\begin{bmatrix} \phantom{i}n\phantom{i}\\ \phantom{i}k\phantom{i} \end{bmatrix} \smashoperator[r]{\sum_{\sum_1^k m_i=n-1}}c_{1+m_1}...c_{1+m_k}(n-m_1)...(n-m_1-...-m_k)
$$
It is left to notice that in addition to $\psi(x)=\mathrm{Q}\mathfrak{T}^{-1}xe^{-x}$, we have $\psi(x)=x\frac{\partial}{\partial x}\mathrm{Q}\mathfrak{T}^{-1} \ln(1+x)$.

\begin{center}
\textbf{Conclusion}
\end{center}
The main goal of this study was to investigate the equation $f/f'=f(px)/p$. The fact, that solutions of this equation in case $p^4=\pm 4$ are actually known, came as a surprise to me, but now, when the answer to this question is known, it doesn't seem really that surprising. Unfortunately, we are far from full understanding of the structure of this equation, but it is hoped, that one can say something else about the resulting series apart from their coefficients.

\newpage

\begin{center}
\textbf{Appendix A}
\end{center}

\hrule

~\\
~\\
Suppose $\gamma_p=(\Delta_p e^{-x})^{inv}$. Then the following holds:\\
\strut\dotfill\\
$$
\begin{CD}
\displaystyle\int_{0}^{x}\frac{dt}{\sqrt[n]{1-t^n}} @>\mathrm{Q}\mathfrak{T}\mathrm{Q}>>\displaystyle\int_{0}^{x}\frac{dt}{1+t^n} @>\mathrm{Q}\mathfrak{T}\mathrm{Q}>> \displaystyle\int_{0}^{x}\gamma'_{\frac{1}{n}}(nt^n)e^{\frac{1-n}{n}\gamma_{\frac{1}{n}}(nt^n)}dt@>\mathrm{Q}\mathfrak{T}\mathrm{Q}>> 
\end{CD}
$$
$$
\begin{CD}
@>\mathrm{Q}\mathfrak{T}\mathrm{Q}>> \displaystyle\int_{0}^{x} \frac{(1+n(1-n)\gamma'_{\nicefrac{1}{n}}(n(1-n)t^n))^2}{1+n(2-n)\gamma'_{\nicefrac{1}{n}}(n(1-n)t^n)}dt
\end{CD}
$$
\strut\dotfill\\
$$
\begin{CD}
\displaystyle\int_{0}^{x}\frac{dt}{(1-t^n)^{\frac{2}{n}}} @>\mathrm{Q}\mathfrak{T}\mathrm{Q}>>\displaystyle\int_{0}^{x}\frac{dt}{\sqrt{1+4t^n}} @>\mathrm{Q}\mathfrak{T}\mathrm{Q}>> \displaystyle\int_{0}^{x}\gamma'_{\frac{2}{n}}(2nt^n)e^{\frac{2-n}{n}\gamma_{\frac{2}{n}}(2nt^n)}dt
\end{CD}
$$
\strut\dotfill\\
$$
\begin{CD}
\displaystyle\int_{0}^{x}e^{\gamma_p(t)}dt @>\mathrm{Q}\mathfrak{T}\mathrm{Q}>>\displaystyle\int_{0}^{x}\frac{1+(p-1)t}{1+pt}dt @>\mathrm{Q}\mathfrak{T}\mathrm{Q}>>
\end{CD}
$$
$$
\begin{CD}
@>\mathrm{Q}\mathfrak{T}\mathrm{Q}>> \displaystyle\int_{0}^{x}\frac{1+pt+\sqrt{1+2(p-2)t+p^2t^2}}{2\sqrt{1+2(p-2)t+p^2t^2}}dt
\end{CD}
$$
\strut\dotfill\\
$$
\begin{CD}
\displaystyle\int_{0}^{x}(1+t)^{-\frac{1}{p}}(1+(1-p)t)^{\frac{1}{p}-1}dt @>\mathrm{Q}\mathfrak{T}\mathrm{Q}>>\displaystyle\int_{0}^{x}\gamma_p'^2(t) e^{(p-2)\gamma_p(t)}dt
\end{CD}
$$
\strut\dotfill\\
$$
\begin{CD}
\displaystyle\int_{0}^{x}\gamma_p'^{1-p}(t)e^{p(1-p)\gamma_p(t)}dt @>\mathrm{Q}\mathfrak{T}\mathrm{Q}>>\displaystyle\int_{0}^{x}\gamma_{\frac{1}{p}}'^2(pt) e^{-2\gamma_{\frac{1}{p}}(pt)}dt
\end{CD}
$$
\strut\dotfill\\
$$
\begin{CD}
\displaystyle\int_{0}^{x}\gamma_p'^{1-\frac{1}{\alpha}}(t)e^{\frac{1-p}{\alpha}\gamma_p(t)}dt @>\mathrm{Q}\mathfrak{T}\mathrm{Q}>>\displaystyle\int_{0}^{x}\frac{\gamma_{\alpha}'( \frac{t}{\alpha})}{1+(1-p)te^{\gamma_{\alpha}(\frac{t}{\alpha})}} e^{-\gamma_{\alpha}(\frac{t}{\alpha})}dt
\end{CD}
$$
\strut\dotfill\\
$$
\begin{CD}
\displaystyle\int_{0}^{x}\gamma_p'^{-1}(t)e^{(2-p)\gamma_p(t)}dt @>\mathrm{Q}\mathfrak{T}\mathrm{Q}>>\displaystyle\int_{0}^{x}\frac{1+(p-2)t+\sqrt{1+2(p-2)t+p^2t^2}}{2\sqrt{1+2(p-2)t+p^2t^2}}dt
\end{CD}
$$
\strut\dotfill\\
$$
\begin{CD}
\displaystyle\int_{0}^{x}e^{\alpha\gamma_p(t)}dt @>\mathrm{Q}\mathfrak{T}\mathrm{Q}>>\displaystyle\int_{0}^{x}\gamma_{\tfrac{p}{1-\alpha}}((1-\alpha)t)\left[(p-1)t+e^{-\gamma_{\frac{p}{1-\alpha}}((1-\alpha)t)}\right]dt
\end{CD}
$$

\newpage

\begin{center}
\textbf{Appendix B}
\end{center}

\hrule

~\\
~\\
\underline{\textbf{Theorem}}. Consider $f(x)=\Delta_p(1+A\Delta_{-p})$. Suppose $A_f=\alpha f'(D)^{-1}$, $D_f=f(D)$. Then if we denote the falling factorials by $(x)_n$, for any integral $n$ the following holds:
$$
A_f^{n+1} f'(D)^{2n+1} A_f^{n\vphantom{+1}}=p^{2n+1}\left(\frac{\alpha}{p}+n\right)_{2n+1}
$$
\emph{Proof}: It is sufficient to prove the statement for any $p\neq 0$, since the case $p=0$ is the degeneration of general one. Expand the brackets and put $\tau=A/p$. Then $f'(x)=$ $=(1+\tau)e^{px}-\tau e^{-px}$, and hence $(f'(D))^2=1+2p(1+2\tau)D_f+p^2 D_f^2$, $f''=p^2 f+p(1+2\tau)$. It is sufficient to show that if $g_n=A_f^{n+1} f'(D)^{2n+1} A_f^{n\vphantom{+1}}$, then $g_{n\vphantom{-1}}^{\vphantom{-1}} g_{n-1}^{-1}=\alpha^2-n^2p^2$, since $g_0=\alpha$.
\begin{align*}
g_{n\vphantom{-1}}^{\vphantom{-1}} g_{n-1}^{-1}&=A_f^{n+1}f'(D)^{2n+1} A_f^{n\vphantom{+1}}A_f^{1-n}f'(D)^{1-2n} A_f^{-n\vphantom{+1}}=\\
&=A_f^{n+1}f'(D)^{2n+1} \alpha f'(D)^{-2n} A_f^{-n\vphantom{+1}}=A_f^{n+1}((2n+1)f''(D)+\alpha f'(D))A_f^{-n\vphantom{+1}}=\\
&=(2n+1)pA_f^{n+1}(1+2\tau+pD_f^{\vphantom{n}})A_f^{-n\vphantom{+1}}+A_f^{n+2}f'(D)^2A_f^{-n\vphantom{+1}}=\\
&=(2n+1)pA_f^{n+1}(1+2\tau+pD_f^{\vphantom{n}})A_f^{-n\vphantom{+1}}+A_f^{n+2}(1+2p(1+2\tau)D_f^{\vphantom{n}}+p^2 D_f^{2\vphantom{n}})A_f^{-n\vphantom{+1}}
\end{align*}
Now use the equality of operators $D_f^{\vphantom{n}} A_f^n=nA_f^{n-1}+A_f^n D_f^{\vphantom{n}}$, to rewrite the latter expression in the following manner:
\begin{align*}
&(1+2\tau)(1+2n)pA_f-n(2n+1)p^2+p^2(2n+1)A_fD_f+A_f^2-\\
&-2np(1+2\tau)A_f+2p(1+2\tau)A_f^2D_f+n(n+1)p^2-2np^2A_fD_f+\\
&+p^2A_f^2D_f^2=\\
&=-n^2p^2+\underline{A_f^2}+p^2A_fD_f+p(1+2\tau)A_f+\underline{2p(1+2\tau)A_f^2D_f}+\underline{p^2A_f^2D_f^2}=\\
&=-n^2p^2+p^2A_fD_f+p(1+2\tau)A_f+A_f^2f'(D)^2
\end{align*}
It is left to notice the explicit formula
$$
A_f^2=\alpha^2\frac{1}{f'(D)^2}-\alpha\frac{f''(D)}{f'(D)^3}=\alpha^2\frac{1}{f'(D)^2}-p\alpha\frac{1+2\tau+pf(D)}{f'(D)^3}
$$
to finally obtain
$$
\boxed{g_{n\vphantom{-1}}^{\vphantom{-1}} g_{n-1}^{-1}=\alpha^2-n^2p^2}
$$
\hfill\ensuremath{\blacksquare}

\end{document}